\documentclass[letterpaper,oneside]{amsart}
\pdfoutput=1 
\usepackage[T1]{fontenc}
\usepackage{color,soul}
\usepackage{amsrefs}
\usepackage[pagewise]{lineno}
\usepackage[dvipsnames]{xcolor}
\usepackage{todonotes}
\setuptodonotes{inline}
\usepackage[utf8]{inputenc} 
\usepackage{comment}
\usepackage{mathtools} 
\usepackage{colonequals}
\usepackage{amsmath}
\usepackage{amssymb}
\usepackage{tikz}
\usetikzlibrary{patterns.meta}
\usepackage{mathrsfs}
\usepackage{enumitem}
\usepackage[foot]{amsaddr}
\usepackage{mathtools}

\usepackage{mleftright}
\usepackage{graphicx}
\usepackage{hyperref}
\usepackage{cancel}
\usepackage[english]{babel}
\usepackage[refpage, notocbasic]{nomencl}
\makenomenclature

\usepackage{tikz}
\usepackage{thmtools}
\usepackage{amssymb}
\usepackage[all]{xy} 

\newcommand{\Id}{\operatorname{Id}}

\newcommand{\Z}{\mathbb{Z}}

\newcommand{\N}{\mathbb{N}}

\newcommand{\Absurd}{(\rightarrow\leftarrow)}

\newcommand{\cat}[1]{\mathbf{#1}}
\newcommand{\ocat}[1]{|\mathbf{#1}|}
\newcommand{\mcat}[3]{[#1,#2]_\mathbf{#3}}

\title[Semi-coarse Spaces: Fundamental Groupoid and van Kampen Theorem]{Semi-coarse Spaces: Fundamental Groupoid and the van Kampen Theorem}
\author{Jonathan Treviño-Marroquín}

\thanks{This work was supported
by the grant N62909-19-1-2134 from the US Office of Naval Research Global and the Southern Office of Aerospace Research and Development
of the US Air Force Office of Scientific Research. The author was also supported by the
CONACYT postgraduate studies scholarship number 839062, and this document is part of his PhD project supervised by Antonio Rieser.}

\theoremstyle{plain}
\newtheorem{teorema}{Theorem}[section]
\newtheorem{proposicion}[teorema]{Proposition}
\newtheorem{corolario}[teorema]{Corollary}
\newtheorem{lema}[teorema]{Lemma}

\theoremstyle{definition}
\newtheorem{definicion}[teorema]{Definition}

\newtheorem{observacion}[teorema]{Remark}
\newtheorem{ejemplo}{Example}[teorema]
\newtheorem{noejemplo}{Non-example}[teorema]

\declaretheoremstyle[
qed=\qedsymbol
]{mystyle}

\begin{document}

\begin{abstract}
In algebraic topology, the fundamental groupoid is a classical homotopy invariant which is defined using continuous maps from the closed interval to a topological space. In this paper, we construct a semi-coarse version of this invariant, using as paths a finite sequences of maps from $\mathbb{Z}_1$ to a semi-coarse space, connecting their tails through semi-coarse homotopy. In contrast to semi-coarse homotopy groups, this groupoid is not necessarily trivial for coarse spaces, and, unlike coarse homotopy, it is well-defined for general semi-coarse spaces. In addition, we show that the semi-coarse fundamental groupoid which we introduce admits a version of the Van Kampen Theorem.
\end{abstract}

\maketitle

\tableofcontents 



\section*{Introduction}
\pagenumbering{arabic}\label{cap.introduccion}

The category of \emph{semi-coarse spaces} is one of several
topological constructs that has been taken up to study algebraic-topological invariants on graphs and semipseudometric spaces, particularly those induced by a metric space and a privileged scale, which together serve as a foundation for topological data analysis. Other
categories which have been proposed for this purpose include \v{C}ech closure spaces
\cites{Bubenik_Milicevic_2021b, Rieser_2021} and Choquet Spaces
\cite{Rieser_arXiv_2022}, which are generalizations topological spaces. Semi-coarse spaces, on the other hand, generalize Roe's axiomatization of coarse
spaces \cite{Roe_2003} by removing the product of sets axiom. This is fundamental because, unlike in coarse spaces, compact (metric) spaces
	admit semi-coarse structures which are no longer necessarily
	equivalent to the point, allowing one to coarsen a compact metric space
	up to a pre-defined scale \cite{Zava_2019} without making it trivial as a semi-coarse space. As a result, the integers provide us with a collection of interval objects for semi-coarse spaces, which we can use to define homotopy in this category \cite{rieser2023semicoarse}. Canonical examples of semi-coarse
spaces are the (undirected) graphs, such as the cyclic graphs and the
Cayley graph of the integers with the unit as generator
($\mathbb{Z}_1$), and semipseudometric spaces. Unfortunately, the
loss of the product axiom for coarse spaces produces two possible
definitions of bounded sets and fails to fulfill that the union of
bounded sets is bounded in connected spaces. However, we may instead use bornological maps as our morphisms, which avoids the need to define bounded sets, but complicates the adaptation of existing notions of coarse homotopy \cite{Mitchener_2020}
to semi-coarse spaces (for more on this, see
\autoref{Appex:CoarseHomotopy}).

In previous work, joint with Antonio Rieser
\cite{rieser2023semicoarse}, we defined homotopy groups for arbitrary semi-coarse
spaces. There we observe that, when applied to coarse spaces, $\pi_0$ quantifies the
number of (coarsely) connected components, and the rest of the homotopy groups
vanish. In this paper, 
we construct the semi-coarse fundamental groupoid, inspired by the one in
topology \cite{Brown_Higgins_Sivera_2011} and by the visual boundary
in semi-geodesic spaces \cite{Kitzmiller_2009}, and show that it may be non-trivial on both coarse and non-coarse semi-coarse spaces. This construction generalizes the
semi-coarse fundamental group, and reduces to it under some reasonable
assumptions. While the definition of the new fundamental groupoid is
slightly technical, it is indeed a semi-coarse invariant which is
non-trivial on coarse spaces and is able to distinguish
them. Moreover, under some specific cases, we are able to prove an analogue to the van
Kampen theorem. In related work, the fundamental groupoid has been also realized for graphs in $A$-theory \cite{kapulkin2023fundamental} and for $\times$-homotopy \cite{Chih_Scull_2022}.

The article is organized as follows. We introduce the category of semi-coarse spaces in \textbf{\autoref{cap.SemiCoarseHomotopy}}, where we also recall some of the results and definitions of \cite{rieser2023semicoarse}, and where we introduce semi-coarse homotopy and the semi-coarse fundamental group. In semi-coarse spaces, we do not have a canonical collection of sets with which to construct useful covers, like open sets in topological spaces, and in \textbf{\autoref{cap.Connectedness_and_Well-Splitting}} we partially resolve this problem by introducing the notion of the well-splitting of a space by a pair of subsets 
(\autoref{def:WellSplit}), which will later provide us with a sufficient condition for a van Kampen theorem for the semi-coarse fundamental groupoid to hold. We then establish several relations between well-splitting, coarse spaces, and connectedness. In 
\textbf{\autoref{cap.SCStrings}}, we define strings of paths and an equivalence relation on them, and we use these objects to construct the semi-coarse fundamental groupoid in \textbf{\autoref{cap.FundGrupoid}}.  

	In \textbf{\autoref{cap.RelativeFundamentalGroupoid}}, we begin with the main results of the paper, and define the relative 
fundamental groupoid and prove a semi-coarse version of the Van Kampen Theorem give some mild hypotheses. Finally, in \textbf{\autoref{Appex:CoarseHomotopy}} we discuss several issues which arise when trying to generalize coarse homotopy \cite{Mitchener_2020} to both coarse and semi-coarse spaces. We do not discard the possibility that this notion can be extended, but it appears that such a generalization is difficult to construct.

\section{Semi-coarse Spaces and its Homotopy}
\label{cap.SemiCoarseHomotopy}

Semi-coarse spaces \cite{Zava_2019} are a generalization of coarse spaces in which ideas of coarse geometry may also be used to study finite undirected graphs and metric spaces with privileged scale. This category was first introduced in \cite{Zava_2019}, and in \cite{rieser2023semicoarse}, the present author, joint with Antonio Rieser, developed the homotopy semi-coarse spaces, adapting ideas of \cite{Babson_etal_2006} to the cartesian product instead the inductive (or box) product as in $A$-theory \cite{Babson_etal_2006}. (Also see \cite{Dochtermann_2009} for a similar construction using the cartesian product, this time on graphs which are not necessarily reflexive.) In this section, we give the basic definitions and summarize the results of \cite{rieser2023semicoarse} which we use in the subsequent sections.

\begin{definicion}[Semi-coarse space; \cite{rieser2023semicoarse}]
\label{def:Semi-coarse} \index{Semi-coarse}
Let $X$ be a set, and let $\mathcal{V}\subset\mathcal{P}(X\times X)$
be a collection of subsets of $X \times X$ which satisfies 
\begin{enumerate}[label=(sc\arabic*)]
    \item \label{sc1}$\Delta_X \in \mathcal{V}$,
    \item \label{sc2} If $B\in\mathcal{V}$ and $A\subset B$, then $A\in\mathcal{V}$,
    \item \label{sc3} If $A,B\in \mathcal{V}$, then $A\cup B\in\mathcal{V}$,
    \item \label{sc4} If $A\in\mathcal{V}$, then $A^{-1}\in\mathcal{V}$.
\end{enumerate}
We call $\mathcal{V}$ a \emph{semi-coarse structure on $X$}, and we say
that the pair $(X,\mathcal{V})$ is a \emph{semi-coarse space}. 
If, in addition, $\mathcal{V}$ satisfies
\begin{enumerate}[label=(sc\arabic*),resume]
    \item\label{item:Coarse condition} If $A,B\in\mathcal{V}$, then $A\circ B\in\mathcal{V}$.
\end{enumerate}
then $\mathcal{V}$ will be called a \emph{coarse structure}, and
$(X,\mathcal{V})$ will be called a \emph{coarse space}, as in \cite{Roe_2003}. 
\end{definicion}

A map $f:X\rightarrow Y$ between the coarse spaces $(X,\mathcal{V})$ and $(Y,\mathcal{W})$ is called \emph{$(\mathcal{V},\mathcal{W})$-bornologous} if $f\times f(V)\coloneqq \{ (f(a),f(b))\mid (a,b)\in V\}\in\mathcal{W}$ for every $V\in\mathcal{V}$. We observe that semi-coarse spaces with bornologous maps is a category. In the consecutive, we will denote this category by $\cat{SCoar}$, the collection of its objects by $\ocat{SCoar}$ and the collection of morphisms from $X$ into $Y$ by $\mcat{X}{Y}{SCoar}$ for every $X,Y\in \ocat{SCoar}$.

Semi-coarse spaces form a topological construct (\cite{Preuss_2002}, 1.1.2), the practical effect of which for this paper is that final and initial structures are defined for every collection of semi-coarse spaces \cite{Zava_2019}. Some important structures initial and final structures are: subspace, product space, quotients and disjoint union space, which we discuss in the results below.

\begin{definicion}[Semi-coarse Subspace; 2.2.3 \cite{rieser2023semicoarse}]
\label{def:SCSubspace}
Let $(X,\mathcal{V})$ be a semi-coarse space and let $Y\subset X$. The pair $(Y,\mathcal{V}_Y)$ is called a \emph{semi-coarse subspace of $X$}, where $\mathcal{V}_Y$ is the collection
\begin{align*}
\mathcal{V}_Y \coloneqq \{ V\cap (Y\times Y) \mid V\in \mathcal{V} \}
\end{align*}
\end{definicion}

In the following sense, bornologous maps may be identified locally: If all of the controlled sets of a semi-coarse space may be written as the finite union of the controlled sets in subspaces $(X_i,\mathcal{V}_i) \subset (X,\mathcal{V})$, then we can observe that if a map is bornologous on each $(X_i,\mathcal{V}_i)$ then they are bornologous in the whole space $(X,\mathcal{V})$. This is formalized in the following.

\begin{proposicion}[\cite{rieser2023semicoarse} 2.2.4]
\label{prop:BornologousPartition}
Let $(X,\mathcal{V})$ and $(Y,\mathcal{W})$ be semi-coarse spaces, and suppose that $(X_i,\mathcal{V}_i)\subset (X,\mathcal{V})$, $i\in\{1,\ldots,n\}$, are subspaces of $(X,\mathcal{V})$ such that $\cup_{i=1}^n X_i = X$ and every set $V\in\mathcal{V}$ may be written in the form
\begin{align}
V = \bigcup_{i=1}^n V_i
\end{align}
where each $V_i\in\mathcal{V}_i$. Now suppose that $f:X\rightarrow Y$ is a map such that the restrictions $f\mid_{X_i}:(X_i,\mathcal{V}_i)\rightarrow (Y,\mathcal{W})$ are bornologous for all $i\in\{1,\ldots,n\}$. Then $f:(X,\mathcal{V})\rightarrow (Y,\mathcal{W})$.
\end{proposicion}

We now define the semi-coarse product, quotient, and disjoint union.

\begin{definicion}[Semi-coarse Product; 2.3.7 \cite{rieser2023semicoarse}]
\label{def:SCProduct}
Let $(X,\mathcal{V})$ and $(Y,\mathcal{W})$ be semi-coarse spaces. $(X\times Y, \mathcal{V}\times \mathcal{W})$ is called \emph{the product of $(X,\mathcal{V})$ and $(Y,\mathcal{W})$}, and $\mathcal{V}\times \mathcal{W}$ \emph{the product structure on $X \times Y$}, where $U\in \mathcal{V}\times \mathcal{W}$ if there exists $V\in\mathcal{V}$ and $W\in\mathcal{W}$ such that $U\subset V\boxtimes W$ with
\begin{align*}
V\boxtimes W \coloneqq \{ ((v,w),(v',w')) \mid (v,v')\in V, (w,w')\in W \}
\end{align*}
\end{definicion}

Since semi-coarse spaces are a topological construct, the inductive product is also defined and it defines another kind of homotopy (see examples of homotopies induced by the inductive product in other categories in \cites{Babson_etal_2006, Bubenik_Milicevic_2021b}.) For the objectives of the current paper, it is enough to work with (categorical) product space. It could be interesting to explore the homotopy theories on semi-coarse spaces induced by the inductive product  in future work.

\begin{definicion}[Quotient Space; 2.4.1 and 2.4.2 \cite{rieser2023semicoarse}]
\label{def:SCQuotient}
Let $(X,\mathcal{V})$ be a semi-coarse space, let $Y$ be a set, and let $g:X\rightarrow Y$ be a surjective function. We define
\begin{align*}
\mathcal{V}_g \coloneqq \{ (g\times g)(V)\mid V\in \mathcal{V}\}.
\end{align*}
Then $(Y,\mathcal{V}_g)$ is called the \emph{quotient space} (or the \emph{semi-coarse space inductively generated by the function $g$}) and $\mathcal{V}_g$ is called the \emph{quotient structure} (or the \emph{semi-coarse structure inductively generated by the function $g$}).
\end{definicion}

\begin{teorema}[\cite{rieser2023semicoarse} 2.4.3]
Let $(X,\mathcal{V})$ be a semi-coarse space, let $(Y,\mathcal{V}_g)$ be the semi-coarse quotient space inductively generated by the function $g:X \rightarrow Y$, and let $(Z,\mathcal{Z})$ be a semi-coarse space. A function $f:(Y,\mathcal{V}_g)\rightarrow (Z,\mathcal{Z})$ is bornologous if and only if $f\circ g: (X,\mathcal{V})\rightarrow (Z,\mathcal{Z})$ is bornologous.
\end{teorema}

\begin{definicion}[Disjoint Union; 2.5.1 and 2.5.2 \cite{rieser2023semicoarse}]
\label{def:SCDisjointUnion}
Let $\{ (X_\lambda, \mathcal{V}_\lambda) \}_{\lambda\in \Lambda}$ be a collection of semi-coarse spaces indexed by the set $\Lambda$, and let $\sqcup_{\lambda\in \Lambda} \mathcal{V}_\lambda$ be the collection of sets of the form
\begin{align*}
\bigsqcup_{\lambda \in \Lambda'} A_\lambda,
\end{align*}
where $|\Lambda'|<\infty$ and each $A_\lambda\in \mathcal{V}_\lambda$. (Note that any given $A_\lambda$ may be the empty set.)

We call $\sqcup_{\lambda\in \Lambda} \mathcal{V}_\lambda$ the disjoint union semi-coarse structure, and $( \sqcup_{\lambda\in \Lambda} X_\lambda , \sqcup_{\lambda\in \Lambda} \mathcal{V}_\lambda )$ the disjoint union of the semi-coarse spaces $\{ (X_\lambda , \mathcal{V}_\lambda) \}_{\lambda\in \Lambda}$.
\end{definicion}

\begin{proposicion}
\label{prop:DisjointUnion}
Let $\{(X_\lambda,\mathcal{V}_\lambda\}$ be a collection of semi-coarse spaces indexed bu the set $\Lambda$, and let $(Z,\mathcal{Z})$ be a semi-coarse space. Define $i_\lambda: X_\lambda\rightarrow X$ such that $i_\lambda(x)=(x,\lambda)$. Then $g:(\sqcup X_\lambda, \sqcup \mathcal{V}_\lambda)\rightarrow (Z,\mathcal{Z})$ is a bornologous maps if and only if $g\circ i _\lambda: X_\lambda\rightarrow Z$ is bornologous for every $\lambda\in\Lambda$.
\end{proposicion}

\begin{proof}
By definition, $i_\lambda$ sends controlled sets to controlled sets, so it is a bornologous map. Therefore, if $g$ is bornologous, we have that $g\circ i_\lambda$ is bornologous for every $\lambda\in \Lambda$.

On the other hand, consider that $g\circ i_\lambda$ is bornologous for every $\lambda$. Take a controlled set $A$ of $\sqcup X$, then there exists $\lambda_1,\ldots,\lambda_n\in \Lambda$ and $A_{\lambda_j}\in \mathcal{V}_j$ such that $A= \sqcup A_{\lambda_j}$. Since $g\circ i_{\lambda_j}$ is bornologous, then $g(A_{\lambda_j})\in\mathcal{Z}$, and therefore $\bigcup g(A_{\lambda_j})\in\mathcal{Z}$. It follows that $g$ is bornologous.
\end{proof}

Through the disjoint union, we can generate the coarser coarse space which contains a semi-coarse space.

\begin{definicion}[Product Extension; 2.5.9 and 2.5.10 \cite{rieser2023semicoarse}]
\label{def:SCVnPE}
Let $(X,\mathcal{V})$ be a semi-coarse space, and define
\begin{align*}
\mathcal{V}^{PE}\coloneqq \{ C\subset X\times X : \exists A,B\in \mathcal{V} \text{ with } C\subset A\circ B \}.
\end{align*}
We call the structure $\mathcal{V}^{PE}$ the \emph{set product extension of $\mathcal{V}$}, and the ordered pair $(X,\mathcal{V}^{PE})$ is called the \emph{set product extension of $(X,\mathcal{V})$}. For any $k\in\mathbb{N}$, we recursively define $\mathcal{V}^{kPE}$ to be the set product extension of $\mathcal{V}^{(k-1)PE}$.
\end{definicion}

Observe that $\mathcal{V}\subset \mathcal{V}^{PE}$.

\begin{definicion}[Product Extension; 2.5.14 and 2.5.15 \cite{rieser2023semicoarse}]
\label{def:SCInfPE}
Let $(X,\mathcal{V})$ be a semi-coarse space, and let $\{(X,\mathcal{V}^{kPE}, \iota_i^j, \mathbb{N}\}$ be the directed system of semi-coarse spaces such that all the $\iota: (X,\mathcal{V}^{iPE})\rightarrow (X,\mathcal{V}^{jPE})$ are the identity map. We call the direct limit coarse structure \emph{the coarse structure induced by $\mathcal{V}$}, which we denote by $\mathcal{V}^\infty$.
\end{definicion}

To close this section, we introduce briefly semi-coarse homotopy and the semi-coarse fundamental group $\pi_1(X,\mathcal{V},x_0)$. Semi-coarse homotopy will be necessary to define what a string is in \autoref{def:String}; we will later observe that the semi-coarse fundamental groupoid will contain a copy of $\pi_1(X,\mathcal{V},x_0)$ for every basepoint $x_0 \in X$ (whether or not $X$ is connected).

In what follows, we denote by $\mathbb{Z}_1$ the semi-coarse space on $\mathbb{Z}$ such that the controlled sets are subsets of $\{(i,i+j)\mid i\in\mathbb{Z}, j\in \{-1,0,1\}\}$.

\begin{definicion}[Semi-coarse Homotopy; 3.1.3  \cite{rieser2023semicoarse}]
\label{def:SCHomotopy}
Let $X$ and $Y$ be semi-coarse spaces. We say that two bornologous maps $f,g:X\rightarrow Y$ are \emph{homotopic} if there exists a bornologous map $H:X\times \mathbb{Z}_1 \rightarrow Y$ and integers $M,N\in\mathbb{Z}$ such that $H(x,k)=f(x)$ for every $k\leq M$ and $H(x,k)=g(x)$ for every $k\geq N$.
\end{definicion}

\begin{definicion}[Homotopy of Maps of Pairs and Triples; 3.1.6 \cite{rieser2023semicoarse}]
Let $X$ and $Y$ be semi-coarse spaces, let $B\subset A\subset X$ and $D\subset C\subset Y$ be endowed with their respective subspace semi-coarse structures, and let $f,g:(X,A,B)\rightarrow (Y,C,D)$ be bornologous maps of a triple, i.e. maps such that $f\mid_A \subset C$ and $f\mid_B \subset D$. We say that \emph{$f$ is relatively homotopic to $g$} and write $f\simeq_{sc} g$ if and only if there is a homotopy $H:X\times \mathbb{Z}_1\rightarrow Y$ such that $H\mid_{A\times \mathbb{Z}} \subset C$ and $H\mid_{B\times \mathbb{Z}} \subset D$.

We define a homotopy between maps of pairs $f,g:(X,A)\rightarrow (Y,C)$ to be the homotopy between maps of a triple as above, where $B=A$ and $C=D$.
\end{definicion}

We now summarize the construction of the fundamental group from \cite{rieser2023semicoarse}, Section 3.2: Fix a semi-coarse space $X$, and, for every $n \in \N$, consider the set $[n]\coloneqq \{0,\ldots,n\}$ as a subspace of $\mathbb{Z}_1$, and consider all bornologous maps $f:[n]\rightarrow X$ such that $f(0)=f(n)=x$ for a fixed $x\in X$. For every $m\in\mathbb{N}_0$ such that $n<m$, we may extend any bornologous map $f:[n]\rightarrow X$ to the bornologous map $f_m:[m]\rightarrow X$, where $f_m(i)= f(i)$ for every $0\leq i \leq n$ and $f_m(i)=f(n)$ for every $n<i\leq m$.

Define the operator $i_n^m:\mcat{[n]}{X}{SCoar} \to \mcat{[m]}{X}{SCoar}$ by $i_{n}^m f = f_m$. Note that $i^m_n$ induces a map on homotopy classes, which, abusing notation, we also write as $i^m_n[f] = [f_m]$. We now have a directed set
\begin{align*}
([f:[n]\rightarrow X], i_{n}^m,\mathbb{N}),
\end{align*}
and we call the direct limit of this system $\pi_1^{sc}(X,x)$. We denote the classes in $\pi_1^{sc}(X,x)$ by $[f]$, and we say that two classes are homotopic $[f]\simeq_{sc} [g]$ if there exists $f'\in [f]$ and $g'\in [g]$ such that $f' \simeq_{sc} g'$.

Finally, we define a product on these functions: Let $m,m'\in\mathbb{N}$, suppose that $f:\{0,\ldots,m\}\rightarrow X$ and $g:\{0,\ldots, m'\}\rightarrow X$, and define the $\star$-product $f\star g: I_{m+m'}\rightarrow X$ by
\begin{align*}
f\star g(j) \coloneqq \left\lbrace\begin{array}{ll}
f(j) & \text{ if } 0\leq j \leq m\\
g(j-m) & \text{ if } j>m.
\end{array}\right.
\end{align*}

This map induces a product in $\pi_1^{sc}(X,x)$ and by Theorem 3.2.18 in \cite{rieser2023semicoarse} we have that $\pi_1^{sc}(X,x)$ is a group with $\star$-operation.

\section{Connectedness and Well-Splitting}
\label{cap.Connectedness_and_Well-Splitting}

Topological spaces have a considerable advantage over other categories: open sets. When we have an open cover $\{U,V\}$, we already get that the pushout of $U\leftarrow U\cap V \rightarrow V$ is exactly $X$ up isomorphisms; furthermore, we would know that $X$ is disconnected if $U\cap V=\varnothing$. In general, semi-coarse spaces do not have a family of sets with such characteristics. To remedy this deficiency, 
in this section we recall the definition for connectedness for semi-coarse spaces from \cite{rieser2023semicoarse}, definition 3.3.1, and we introduce the property that two subsets $A,B \subset X$ \emph{well-split} $X$ in \autoref{def:WellSplit}. As we will see, this notion has enough of the properties of an open cover of two elements of a topological space to make it an important tool in the proof of the semi-coarse Van Kampen theorems .

\begin{definicion}
\label{def:PathConnected}
Let $X$ be a semi-coarse space. We say that $X$ is connected if, for every $x,y\in X$, there exists a subspace $\{0, 1, \ldots, n\}$ of $\mathbb{Z}_1$ and a path $\gamma:\{0, 1, \ldots, n\}\rightarrow X$ such that $\gamma(0)=x$ and $\gamma(n)=y$.
\end{definicion}

\begin{proposicion}[Universal Property of Pushout]
\label{prop:PushoutIsPushout}
Let $A,\ B$ and $C$ be semi-coarse spaces, let $f:C\rightarrow A$ and $g:C\rightarrow B$ be bornologous maps, and let $A\sqcup_C B\coloneqq A\sqcup B/(f(c)\sim g(c))$. Then $A \sqcup_C B$, with the semi-coarse structure induced by the disjoint union and the quotient $\mathcal{V}_{A \sqcup_C B}$, is the pushout in the commutative diagram
\begin{align*}
\xymatrix{
C \ar[r]^{f} \ar[d]_{g} & A \ar[d]^{i_1} \\
B \ar[r]^{i_2} & A\sqcup_C B
}
\end{align*}
i.e., for every other commutative diagram
\begin{align*}
\xymatrix{
C \ar[r]^{f} \ar[d]_{g} & A \ar[d]^{j_1} \\
B \ar[r]^{j_2} & X
}
\end{align*}
we have that there exists a unique bornologous map $h:A\sqcup_C B\rightarrow X$ such that $j_1= h i_1$ and $j_2= h i_2$.
\end{proposicion}

\begin{proof}
We use the identification $A\sqcup B \cong (A\times \{1\})\cup (B\times \{2\})$ to aid in the differentiation of elements of $A$ and $B$ inside $A \sqcup B$. Moreover, we define $i_1:A\to A\sqcup B$ by $i_1(a)=(a,1)$ and $i_2:B \to A \sqcup B$ by $i_2(b)=(b,2)$.

We now define the map $h:A\sqcup B \to X$  by
\begin{equation*}
h(x,i) \coloneqq \begin{cases} j_1(x)& (x,i) \in A\times {1}\\
	j_2(x) & (x,i) \in B \times {2}
	\end{cases} 
\end{equation*}
and we define $\overline{h}: A\sqcup_C B\rightarrow X$ by 
$\overline{h}[x,i] = h(x,i)$. We must now show that $\overline{h}$ is a well-defined bornologous map, and that it is the unique bornologous map such that $j_2 = \overline{h}\circ i_2$ and $j_1 = \overline{h}\circ i_1$.

To prove that $\overline{h}$ is well-defined, we examine the following cases:
\begin{enumerate}
\item If $(a,1)\in [b,2]$, then there exists $c\in C$ such that $f(c)=a$ and $g(c)=b$. Thus $\overline{h}[a,1]= \overline{h}i_1(a)= \overline{h}i_1f(c)= \overline{h}i_2g(c)= \overline{h}i_2(b)= \overline{h}[b,2]$.
\item If $(a',1) \in [a,1]$, then there exists $c,c'\in C$ such that $f(c)= a$, $f(c')= a'$ and $g(c)= g(c')$. Thus $\overline{h}[a,1]= \overline{h}[g(c),2]= \overline{h}[g(c'),2]= \overline{h}[a',1]$.
\end{enumerate}
For the cases $(b,2)\in [a,1]$ and $(b',2)\in [b,2]$, the procedure is completely analogous to the above. Therefore, $\overline{h}$ is 
	well-defined.

We prove that $\overline{h}$ is bornologous using the universal properties of the disjoint union and quotient, and observing that the following diagrams are commutative:¨
\begin{align*}
\xymatrix{
A \ar[r]^{i_1} \ar[rd]^{j_1} & A\sqcup B \ar[d]^{h} & B \ar[l]_{i_2} \ar[dl]_{j_2} \\
& X
}\
\xymatrix{
A\sqcup B \ar[r]^q \ar[d]^{h} &  A\sqcup_C B \ar[ld]^{\overline{h}}\\
X
}
\end{align*}
Finally, we observe that $\overline{h}$ is the unique map which satisfies $j_2 = \overline{h}\circ i_2$ and $j_1 = \overline{h}\circ i_1$: Suppose, on the contrary, that there exists $h':A\sqcup_C B\rightarrow X$ such that $h'[x,k]\neq \overline{h}[x,k]$, for some $x\in A \sqcup B$. Then $h'i_k(x)\neq \overline{h}[i_k(x)]=h\circ i_k(x) = j_k(x)$.
\end{proof}

\begin{proposicion}
\label{prop:SplittingWell_WhichSide}
Let $(X,\mathcal{V})$ a semi-coarse space. If $\{(x_1,x_2)\}\in\mathcal{V}_A \sqcup_{\mathcal{V}_{A\cap B}}\mathcal{V}_B$ and:
\begin{enumerate}
\item $x_1\in A$, $x_2\notin A$, then $x_1\in A \cap B$.
\item $x_1\in B$, $x_2\notin B$, then $x_1\in A \cap B$.
\end{enumerate}
\end{proposicion}

\begin{proof}
Let $(X,\mathcal{V})$ a semi-coarse space, $A,B\subset X$, and $\{(x_1,x_2)\}\in\mathcal{V}_{A\sqcup_{A\cap B} B}$. Both cases above are completely analogous; thus we prove just the first case. Consider $x_1,x_2 \in X$ with $x_1\in A$ and $x_2\notin A$.

Let $(Y,\mathcal{W})$ be the subspace $\{1,2,3\}$ of $\mathbb{Z}_1$, and define the maps $j_1:A\rightarrow Y$ and $j_2:B \rightarrow Y$ by
\begin{align*}
j_1(a) &\coloneqq \begin{cases} 2 & a \in A\cap B,\\
	1 & a \in A-B,
\end{cases}\\
j_2(b) &\coloneqq \begin{cases} 2 & b \in A\cap B,\\
	3 & b \in B-A.
	\end{cases}
\end{align*}
Since $\{1,2\}\times \{1,2\}, \{2,3\}\times \{2,3\}\in \mathcal{W}$, then $j_1,j_2$ are bornologous maps. Moreover, the diagram
\begin{align*}
\xymatrix{
A\cap B \ar[r]^{i_1} \ar[d]_{i_2} & A \ar[d]^{j_1} \\
B \ar[r]^{j_2} & Y
}
\end{align*}
is commutative. Thus, since $A\sqcup_{A\cap B} B$ is a pushout, there exists a unique $h: A\sqcup_{A\cap B} B\rightarrow Y$ bornologous such that $j_2=h i_2$ and $j_1 = hi_1$. Since $h$ is bornologous, then $\{h(x_1),h(x_2)\}\in\mathcal{W}$. In addition, $x_1\in A$ and $x_2\in B-A$ by hypothesis, thus $h(x_2)= 3$ and $h(x_1)\in \{1,2\}$. Therefore, $h(x_1)\neq 1$ because $\{(1,3)\}\notin \mathcal{W}$, and it follows that $h(x_1)=2$, so $x_1\in A\cap B$.
\end{proof}

\begin{corolario}
\label{prop:PathBetweenElements}
Let $X$ be a semi-coarse space and $A, B \subset X$. Then every path $\gamma$ whose image is in $A\sqcup_{A\cap B} B$ satisfies that
\begin{itemize}
\item if $\gamma(i)\in A$ and $\gamma(i+1)\notin A$, then $\gamma(i)\in A\cap B$, and
\item if $\gamma(i)\in B$ and $\gamma(i+1)\notin B$, then $\gamma(i)\in A\cap B$.
\end{itemize}
\end{corolario}

\begin{definicion}
\label{def:WellSplit}
Let $(X,\mathcal{V})$ be a semi-coarse space and $A,B\subset X$ non-empty with $X = A \cup B$. Then $A,B$ well-split $X$ if for every $x,y,y'\in X$ such $\{(x,y),(x,y')\}\in\mathcal{V}-\mathcal{V}_{A\sqcup_{A\cap B} B}$ and $\{(y,y')\}\in \mathcal{V}_{A\sqcup_{A\cap B} B}$:
\begin{enumerate}
\item There exist two bornologous paths $\gamma,\gamma':\{0, 1, 2\}\rightarrow (A\sqcup_{A \cap B} B,\mathcal{V}_{A \sqcup_{A\cap B} B})$, where $\{0,1,2\} \subset \Z_1$, such that $\gamma(0)=\gamma'(0)=x$, $\gamma(2)=y$, $\gamma'(2)=y'$ and $\gamma(1)=\gamma'(1)$.

\item The set $\{ \gamma(1)\mid \gamma \text{ is a path in } A\sqcup_{A\cap B} B \text{ such that} \{ \gamma(0)=x,\gamma(2)=y \}$ is connected as subspace of $X$.
\end{enumerate}

\end{definicion}

Condition two is equivalent to saying that every $f:\{0,1,2\}\rightarrow X$ with $f(0)=x$ and $f(2)=y$ are homotopic. This implies that, in a coarse space $X$, every $A,B$ which cover $X$ and which satisfy the first condition in well-split the space. Indeed, the reason to include both of these conditions is to include coarse spaces in our analysis.

\begin{ejemplo}
Let $X\coloneqq \{ x,x',y,y',w,w'\}$ with the semi-coarse structure induced by the following the graph:

\begin{center}
\begin{tikzpicture}
\coordinate (x1) at (0,0);
\coordinate (x2) at (-1,0);
\coordinate (y1) at (1,1);
\coordinate (y2) at (0,1);
\coordinate (w1) at (1,0);
\coordinate (w2) at (-1,1);

\filldraw[black] (x1) circle (0pt) node[anchor=north]{$x$};
\filldraw[black] (x2) circle (0pt) node[anchor=north]{$x'$};
\filldraw[black] (y1) circle (0pt) node[anchor=south]{$y$};
\filldraw[black] (y2) circle (0pt) node[anchor=south]{$y'$};
\filldraw[black] (w1) circle (0pt) node[anchor=north]{$w$};
\filldraw[black] (w2) circle (0pt) node[anchor=south]{$w'$};

\draw[line width=0.5mm,black] (w1)--(x1);
\draw[line width=0.5mm,black] (w1)--(y1);
\draw[line width=0.5mm,black] (w1)--(y2);
\draw[line width=0.5mm,black] (w1)--(w2);

\draw[line width=0.5mm,black] (x1)--(x2);
\draw[line width=0.5mm,Black] (x1)--(y1);
\draw[line width=0.5mm,Black] (x1)--(y2);
\draw[line width=0.5mm,black] (x1)--(w2);

\draw[line width=0.5mm,Black] (x2)--(y2);
\draw[line width=0.5mm,black] (x2)--(w2);

\draw[line width=0.5mm,black] (y2)--(y1);
\draw[line width=0.5mm,black] (y2)--(w2);

\filldraw[black] (x1) circle (2pt);
\filldraw[black] (x2) circle (2pt);
\filldraw[black] (y1) circle (2pt);
\filldraw[black] (y2) circle (2pt);
\filldraw[black] (w1) circle (2pt);
\filldraw[black] (w2) circle (2pt);
\end{tikzpicture}
\end{center}
and let $A\coloneqq \{x,x',w,w'\}$ and $B=\{y,y',w,w'\}$. We observe that we have only three edges outside of $\mathcal{V}_{A \sqcup_{A\cap B} B}$, they are $\{x',y'\}$, $\{x,y'\}$ and $\{x,y\}$, and there are just two pair with a vertex in common, which are: $\{ \{x',y'\} , \{x,y'\} \}$ and $\{ \{x,y'\} , \{x,y\} \}$. Thus, the two conditions are easy to follow and $A$ and $B$ well split $X$. \autoref{fig:PictureSplitingWell} illustrates how $X$ the vertices and edges are divided.
\end{ejemplo}

\begin{figure}[h!]
\begin{tikzpicture}
\coordinate (x1) at (0,0);
\coordinate (x2) at (-1,0);
\coordinate (y1) at (1,1);
\coordinate (y2) at (0,1);
\coordinate (w1) at (1,0);
\coordinate (w2) at (-1,1);

\filldraw[black] (x1) circle (0pt) node[anchor=north]{$x$};
\filldraw[black] (x2) circle (0pt) node[anchor=north]{$x'$};
\filldraw[black] (y1) circle (0pt) node[anchor=south]{$y$};
\filldraw[black] (y2) circle (0pt) node[anchor=south]{$y'$};
\filldraw[black] (w1) circle (0pt) node[anchor=north]{$w$};
\filldraw[black] (w2) circle (0pt) node[anchor=south]{$w'$};

\draw[line width=0.5mm,Red] (w1)--(x1);
\draw[line width=0.5mm,RoyalBlue] (w1)--(y1);
\draw[line width=0.5mm,RoyalBlue] (w1)--(y2);
\draw[line width=0.5mm,Orchid] (w1)--(w2);

\draw[line width=0.5mm,Red] (x1)--(x2);
\draw[line width=0.5mm,Black] (x1)--(y1);
\draw[line width=0.5mm,Black] (x1)--(y2);
\draw[line width=0.5mm,Red] (x1)--(w2);

\draw[line width=0.5mm,Black] (x2)--(y2);
\draw[line width=0.5mm,Red] (x2)--(w2);

\draw[line width=0.5mm,RoyalBlue] (y2)--(y1);
\draw[line width=0.5mm,RoyalBlue] (y2)--(w2);

\filldraw[Red] (x1) circle (2pt);
\filldraw[Red] (x2) circle (2pt);
\filldraw[RoyalBlue] (y1) circle (2pt);
\filldraw[RoyalBlue] (y2) circle (2pt);
\filldraw[Orchid] (w1) circle (2pt);
\filldraw[Orchid] (w2) circle (2pt);
\end{tikzpicture}
\caption{Well-splitting example: Red points are $A-B$, blue points are $B-A$ and violet points are $A\cap B$. $X$ with the semi-coarse space induce by the graph.}
\label{fig:PictureSplitingWell}
\end{figure}

\begin{noejemplo}
Take $X=\mathbb{Z}_1$ and let $A$ be the odd integers and $B$ the even integers. We observe that there is no continuous path from $1$ to $2$ such that every edge is within $\mathcal{V}_{A\sqcup_{A\cap B} B}$. Condition $(2)$ is satisfied by vacuity. \autoref{fig:NonSplitingWell_1} illustrates how $X$ the vertices and edges are divided.
\end{noejemplo}

\begin{figure}[h!]
\begin{tikzpicture}
\coordinate (a') at (-2.5,0);
\coordinate (a) at (-2,0);
\coordinate (b) at (-1,0);
\coordinate (c) at (0,0);
\coordinate (d) at (1,0);
\coordinate (e) at (2,0);
\coordinate (e') at (2.5,0);

\filldraw[black] (a) circle (0pt) node[anchor=north]{$-2$};
\filldraw[black] (b) circle (0pt) node[anchor=north]{$-1$};
\filldraw[black] (c) circle (0pt) node[anchor=north]{$0$};
\filldraw[black] (d) circle (0pt) node[anchor=north]{$1$};
\filldraw[black] (e) circle (0pt) node[anchor=north]{$2$};

\draw[line width=0.5mm,Black] (a')--(a);
\draw[line width=0.5mm,Black] (a)--(b);
\draw[line width=0.5mm,Black] (b)--(c);
\draw[line width=0.5mm,Black] (c)--(d);
\draw[line width=0.5mm,Black] (d)--(e);
\draw[line width=0.5mm,Black] (e)--(e');

\filldraw[RoyalBlue] (a) circle (2pt);
\filldraw[Red] (b) circle (2pt);
\filldraw[RoyalBlue] (c) circle (2pt);
\filldraw[Red] (d) circle (2pt);
\filldraw[RoyalBlue] (e) circle (2pt);

\filldraw[Black] (-2.5-1/6,0) circle (0.2pt);
\filldraw[Black] (-2.5-2/6,0) circle (0.2pt);
\filldraw[Black] (-2.5-3/6,0) circle (0.2pt);

\filldraw[Black] (2.5+1/6,0) circle (0.2pt);
\filldraw[Black] (2.5+2/6,0) circle (0.2pt);
\filldraw[Black] (2.5+3/6,0) circle (0.2pt);
\end{tikzpicture}
\caption{No well-splitting example 1: Red points are $A-B$, blue points are $B-A$ and there are no points in $A\cap B$.}
\label{fig:NonSplitingWell_1}
\end{figure}

\begin{noejemplo}
Let $X=\{a,b,x,y,z\}$ with the semi-coarse structure induced by the following graph:
\begin{center}
\begin{tikzpicture}
\coordinate (a) at (0,1);
\coordinate (b) at (0,-1);
\coordinate (x) at (-1,0);
\coordinate (y) at (1,0);
\coordinate (z) at (2,0);

\filldraw[black] (a) circle (0pt) node[anchor=south]{$a$};
\filldraw[black] (b) circle (0pt) node[anchor=north]{$b$};
\filldraw[black] (x) circle (0pt) node[anchor=east]{$x$};
\filldraw[black] (y) circle (0pt) node[anchor=west]{$y$};
\filldraw[black] (z) circle (0pt) node[anchor=west]{$z$};

\draw[line width=0.5mm,black] (a)--(x);
\draw[line width=0.5mm,black] (a)--(y);
\draw[line width=0.5mm,black] (a)--(z);

\draw[line width=0.5mm,black] (b)--(x);
\draw[line width=0.5mm,black] (b)--(y);
\draw[line width=0.5mm,black] (b)--(z);

\draw[line width=0.5mm,black] (a)--(b);

\filldraw[black] (a) circle (2pt);
\filldraw[black] (b) circle (2pt);
\filldraw[black] (x) circle (2pt);
\filldraw[black] (y) circle (2pt);
\filldraw[black] (z) circle (2pt);
\end{tikzpicture}
\end{center}
and let $A=\{a,x,y,z\}$ and $B=\{b,x,y,z\}$. We can observe that $X$ satisfies the condition $(1)$ and it does not satisfy the condition $(2)$ because $\{x,y,z\}$ are disconnected with semi-coarse subspace structure.
\end{noejemplo}

\begin{figure}[h!]
\begin{tikzpicture}
\coordinate (a) at (0,1);
\coordinate (b) at (0,-1);
\coordinate (x) at (-1,0);
\coordinate (y) at (1,0);
\coordinate (z) at (2,0);

\filldraw[black] (a) circle (0pt) node[anchor=south]{$a$};
\filldraw[black] (b) circle (0pt) node[anchor=north]{$b$};
\filldraw[black] (x) circle (0pt) node[anchor=east]{$x$};
\filldraw[black] (y) circle (0pt) node[anchor=west]{$y$};
\filldraw[black] (z) circle (0pt) node[anchor=west]{$z$};

\draw[line width=0.5mm,Red] (a)--(x);
\draw[line width=0.5mm,Red] (a)--(y);
\draw[line width=0.5mm,Red] (a)--(z);

\draw[line width=0.5mm,RoyalBlue] (b)--(x);
\draw[line width=0.5mm,RoyalBlue] (b)--(y);
\draw[line width=0.5mm,RoyalBlue] (b)--(z);

\draw[line width=0.5mm,Black] (a)--(b);

\filldraw[Red] (a) circle (2pt);
\filldraw[RoyalBlue] (b) circle (2pt);
\filldraw[Orchid] (x) circle (2pt);
\filldraw[Orchid] (y) circle (2pt);
\filldraw[Orchid] (z) circle (2pt);
\end{tikzpicture}
\caption{No well-splitting example 2: Red points are $A-B$, blue points are $B-A$ and violet points are $A\cap B$. $X$ with the semi-coarse space induce by the graph.}
\label{fig:NonSplitingWell_2}
\end{figure}

The next proposition illustrates that well-splitting is a generalization of the condition that $X$ be semi-coarse isomorphic to $A \sqcup_{A \cap B} B$ for some $A,B \subset X$.

\begin{proposicion}
	\label{prop:Homeotowellsplit}
Let $X$ be a semi-coarse space and $A,B\subset X$. If $X$ is (semi-coarse) homeomorphic to $A\sqcup_{A\cap B} B$, then $A,B\subset X$ well-split $X$.
\end{proposicion}

The example in \autoref{fig:PictureSplitingWell} illustrates that the other direction is not necessarily true. Well-splitting will be 
	used as the condition on the cover of $X$ in the semi-coarse van Kampen theorems \ref{theo:VanKampenConnected} and \ref{theo:VankampenDisconnected}. It might be tempting to only use the condition that $A,B\subset X$ cover $X$ and that $X$ and $A\sqcup_{A\cap B}B$ are semi-coarse isomorphic; however, we might leave out the interesting coarse cases. Consider, for example $\mathbb{R}$ with its coarse structure induced by the metric, and $A\coloneqq[0,\infty)$, $B\coloneqq (-\infty,0]$. $A\cup B = X$, but $A\sqcup_{A\cap B} B$ is not a coarse space and thus is not isomorphic to $X$ (see \autoref{fig:PushOut} to compare both spaces). The pair $A,B$ do, however, well-split $X$, since $X$ is coarse, $A\cup B = X$ (so they satisfy condition $(2)$), and they clearly satisfy condition $(1)$ of the definition.

\begin{figure}[h!]
\begin{tikzpicture}[domain=0:4]
	\filldraw[red!24] (0,0)--(0,2)--(2,2)--(2,0)--(2,0)--cycle;
  \filldraw[red!34] (0,0)--(0,1)--(1,2)--(2,2)--(2,1)--(1,0)--cycle;	
	\filldraw[red!49] (0,0)--(0,1.5)--(0.5,2)--(2,2)--(2,0.5)--(1.5,0)--cycle;
  \filldraw[red!70] (0,0)--(0,1)--(1,2)--(2,2)--(2,1)--(1,0)--cycle;
  \filldraw[red] (0,0)--(0,.5)--(1.5,2)--(2,2)--(2,1.5)--(0.5,0)--cycle;
  
  \filldraw[red!24] (0,0)--(0,-2)--(-2,-2)--(-2,0)--(-2,0)--cycle;
  \filldraw[red!34] (0,0)--(0,-1)--(-1,-2)--(-2,-2)--(-2,-1)--(-1,0)--cycle;	
	\filldraw[red!49] (0,0)--(0,-1.5)--(-0.5,-2)--(-2,-2)--(-2,-0.5)--(-1.5,0)--cycle;
  \filldraw[red!70] (0,0)--(0,-1)--(-1,-2)--(-2,-2)--(-2,-1)--(-1,0)--cycle;
  \filldraw[red] (0,0)--(0,-0.5)--(-1.5,-2)--(-2,-2)--(-2,-1.5)--(-0.5,0)--cycle;
  
  \draw[<->] (-2.2,0) -- (2.2,0) node[right] {$\mathbb{R}$};
  \draw[<->] (0,-2.2) node[below]{$A\sqcup_{A\cap B} B$} -- (0,2.2) node[above] {$\mathbb{R}$};
  
  \filldraw[red!8] (4,-2)--(4,2)--(8,2)--(8,-2)--cycle;
  \filldraw[red!12] (4,-2)--(4,1.5)--(4.5,2)--(8,2)--(8,-1.5)--(7.5,-2)--cycle;
  \filldraw[red!17] (4,-2)--(4,1)--(5,2)--(8,2)--(8,-1)--(7,-2)--cycle;
  \filldraw[red!24] (4,-2)--(4,0.5)--(5.5,2)--(8,2)--(8,-0.5)--(6.5,-2)--cycle;
  \filldraw[red!34] (4,-2)--(4,0)--(6,2)--(8,2)--(8,0)--(6,-2)--cycle;
  \filldraw[red!49] (4,-2)--(4,-0.5)--(6.5,2)--(8,2)--(8,0.5)--(5.5,-2)--cycle;
  \filldraw[red!70] (4,-2)--(4,-1)--(7,2)--(8,2)--(8,1)--(5,-2)--cycle;
  \filldraw[red] (4,-2)--(4,-1.5)--(7.5,2)--(8,2)--(8,1.5)--(4.5,-2)--cycle;
  
  \draw[<->] (3.8,0) -- (8.2,0) node[right] {$\mathbb{R}$};
  \draw[<->] (6,-2.2) node[below]{$X$} -- (6,2.2) node[above] {$\mathbb{R}$} ;
\end{tikzpicture}
\caption{Every controlled set in $A\sqcup_{A\cap B}B$ and $X$, respectively, is contained in some red area. Observe that there are not controlled sets in the second and fourth quadrants in $A\sqcup_{A\cap B}B$, in contrast to $X$.}
\label{fig:PushOut}
\end{figure}

The interest in well-splitting is that it gives a set-based condition for
path-connectedness of a semi-coarse space, analogous to that of connectedness in topological spaces. The next proposition makes this clear.

\begin{proposicion}
	\label{prop:SplitWell_IFF_PathConnected}
	Let $(X,\mathcal{V})$ be a semi-coarse space and $A,B\subset X$ non-empty. $X$ is disconnected if and only if there exist $A,B$ which well-split $X$ such that $A\cap B = \varnothing$.
\end{proposicion}

\begin{proof}
	Let $(X,\mathcal{V})$ be a semi-coarse space and $A,B\subset X$ non-empty.
	
	$(\Rightarrow)$ Suppose that $X$ is disconnected and $x\in X$. Take $A$ as the component of $x$ and $B\coloneqq X-A$. Observe that $A\sqcup_{A\cap B} B = A\sqcup B = X$ as semi-coarse spaces. Thus, $A,B$ well-split $X$ by \autoref{prop:Homeotowellsplit}.  
	
	$(\Leftarrow)$ Suppose that $A,B$ well-split $X$ and $A\cap B=\varnothing$. Let $a\in A$ and $b\in B$, and consider that there is a bornologous path $\gamma:\{0,\ldots,n\}\rightarrow X$ such that $\gamma(0)=a$ and $\gamma(n)=b$. Since $A\cup B=X$, there there exists $i$ such that $\gamma(i)\in A$ and $\gamma(i+1)\in B$. However, $\{(\gamma(i),\gamma(i+1))\} \in \mathcal{V}-\mathcal{V}_{A \sqcup_{A\cap B} B}$. Hence, by well-splitting condition $(1)$, there is $\lambda:\{0,1,2\}\rightarrow X$ in $A\sqcup_{A\cap B} B$ such that $\lambda(0)=\gamma(i)$ and $\lambda(2)=\gamma(i+1)$. By \autoref{prop:SplittingWell_WhichSide}, we obtain that at least one of the three elements are in $A\cap B$, which is a contradiction $\Absurd$. Therefore, there is no path between $a,b\in X$; $X$ is disconnected.
\end{proof}

\section{Semi-coarse Strings}
\label{cap.SCStrings}

	Building a groupoid from semi-coarse spaces requires a number of subtle considerations.
	On the one hand, we would like to include the fundamental groups as a part of this groupoid in the usual way, by restricting to just one object and considering its automorphisms.
	On the other hand, our motivation is to find and invariant which doesn't lose all of the structure when we work in coarse spaces, and the fundamental group is trivial in coarse spaces, as shown in 
\cite{rieser2023semicoarse}.
	Our solution will be to build the groupoid out of bornologous maps $f:\mathbb{Z}_1\rightarrow X$ to a semi-coarse space $X$.
	From this point on, we will refer to the restrictions $f\mid_{[n,\infty)}$ and $f\mid_{(-\infty,m]}$ for some $m,n,\in\mathbb{Z}$ as  \emph{tails}, and they will be our main focus in the subsequent construction.
	Occasionally, we refer to $f\mid_{[n,\infty)}$ and $f\mid_{(-\infty,m]}$ as \emph{right tail and left tail}, respectively.
	We will soon see that the tails contain all of the information about objects in our groupoid.
	We begin with the following definitions.
´
\begin{definicion}
\label{def:SymmetricMaps}
	Let $f:\mathbb{Z}_1\rightarrow X$:
\begin{itemize}
	\item $f$ is called symmetric if $f(x)=f(-x)$.
	
	\item $\overline{f}$ is called the \emph{opposite direction} of $f:\mathbb{Z}_1\rightarrow X$ with $\overline{f}(x)=f(-x)$.
	 \nomenclature{$\overline{f}$}{The opposite direction of $f$}
\end{itemize}
\end{definicion}

Through the tails of the symmetric maps, we can develop a equivalence relation.

\begin{definicion}
\label{def:EventuallyEqualMaps}
Let $f,g$ be symmetric maps. $f$ and $g$ are \emph{eventually equal} if there exists $N,M\in\mathbb{N}_0$ such that $f(N+k)=g(M+k)$ for every $k\in\mathbb{N}_0$\nomenclature{$\N_0$}{Non-negative integers}.
\end{definicion}

\begin{proposicion}
\label{prop:SymmIsEqRel}
	The relation of two symmetric maps being eventually equal or homotopic is a equivalence relation.
	We write $\langle f \rangle_\infty$ to refer to the class of the element $f$ with this equivalence class.
	\nomenclature{$\langle f \rangle_\infty$}{The class of maps eventually equal to $f$ with $f$ symmetric}
\end{proposicion}

\begin{observacion}
	Let $f,g:\Z_1\to X$ be symmetric maps. Then, $f\simeq_{sc} g$ if and only if $f\mid_{[0,\infty)}\simeq_{sc} g\mid_{[0,\infty)}$. Thus, the homotopy $H:\Z_1\times \Z_1\to X$ between $f$ and $g$ can be chosen such that $H(-,t):\Z_1\to X$ is a symmetric function.
\end{observacion}

\begin{observacion}
	Let $X$ be a set.
	We recall the indiscrete semi-coarse structure of $X$ defined as the ordered pair $(X,2^X)$, and denoted by $X_\text{ind}$.
	
	For any bornologous map $f:\Z_1\to X_\text{ind}$ we obtain the equality
	\[\langle f \rangle_\infty = \{ g:\Z_1\to X_\text{ind}\}.\]
	That is, every semi-coarse space with the indiscrete structure has just one class of symmetric homotopic maps.
\end{observacion}

	The eventually equal classes will be the objects in the eventual definition of the fundamental groupoid.
	We now define what will be the arrows, or morphisms, between classes of eventually equal symmetric maps.

\begin{definicion}
\label{def:String}
	Let $n\in\mathbb{N}$ and $f_i:\Z \to X$ for $i\in \{0,\ldots, n+1\}$ such that $f_0$ and $f_{n+1}$ are symmetric.
	We say that a $F=(f_1,\ldots,f_n)$ is an \emph{$n$-string from $\langle f_0 \rangle_\infty$ to $\langle f_{n+1} \rangle_\infty$ in $X$} if we have the following:
	For every $i\in \{0,\ldots, n\}$ there exists $N_i,M_i\in\Z$ such that
	\[f_{i,R}(z)\coloneqq f_i(z+N_i),\ f_{i+1,L}(z)\coloneqq f_{i+1}(z+M_i)\]
	satisfies that $f_{i,R}\mid_{[0,\infty)} \simeq_{sc} \overline{f_{i+1,L}}\mid_{[0,\infty)}$.

\emph{Strings} are the elements of the collection of all of the $n$-strings for every $n\in\mathbb{N}$.
\end{definicion}

	This chains are well-defined for every homotopic map $g_0$ to $f_0$ and every homotopic map $g_{n+1}$ to $f_{n+1}$.
	This follows from the condition that $f_{i,R}\mid_{[0,\infty)} \simeq_{sc} \overline{f_{i+1,L}}\mid_{[0,\infty)}$

	We use the following notation in the subsequent sections, especially when we talk about the semi-coarse version of the van Kampen theorem.

\begin{definicion}
\label{def:NotationStrings}\
\begin{itemize}
	\item $\mathcal{F}_{X,A,n,f,g}$ is the collection of all of the $n$-string in $X$ from $\langle f\rangle_\infty$ to $\langle g\rangle_\infty$ such that their tails land in $A\subset X$.
	
	\item If we omit the subset $A$, we consider that $A=X$.
	
	\item If we omit the number $n$ we consider is all of the $n$-strings with $n\in\mathbb{N}$.
	
\item If we omit $f$ and $g$, we consider that the strings can start and end in any element.
\end{itemize}
	When the context is clear, we omit $X$.
\end{definicion}

\begin{observacion}
	Let $F=(f_1,\ldots,f_n)$ be a path from $\langle f \rangle_\infty$ to $\langle g \rangle_\infty$  and $G=(g_1,\ldots,g_m)$ be a path from $\langle g \rangle_\infty$ to $\langle h \rangle_\infty$.
	Then we can observe that there exist $N,M,N',M'\in \Z$ such that
	\[ f_{n,R}(z) \coloneqq f_n(z+N),\ g_L(z) \coloneqq g(z+M),  \]
	\[ g_R(z) \coloneqq g(z+N'), g_{1,L}\coloneqq g_1(z+M'), \]
	and $f_{n,R}\mid_{[0,\infty)} \simeq_{sc} \overline{g_{L}}\mid_{[0,\infty)}$ and $g_{R}\mid_{[0,\infty)} \simeq_{sc} \overline{g_{1,L}}\mid_{[0,\infty)}$.
	
	Without loss of generality, suppose that $M\geq N'$;
	the case when $M<N'$ is completely analogous.
	With the help of the \autoref{tab:StarOperation} and remembering that g is a symmetric function, note that we can define $N''\coloneqq N$,  $M''\coloneqq M+M'-N'$ and the maps
	\[ f_{n,R'}(z)\coloneqq f_n(z+N''),\ g_{1,L'}(z)\coloneqq  g_1(z+M''). \]
	We observe that $f_{n,R'}\mid_{[0,\infty)} \simeq_{sc} g_{1,L'}\mid_{[0,\infty)}$ (see \autoref{tab:StarOperationCon}) and conclude that
	\[(f_1,\ldots,f_n,g_1,\ldots,g_m)\]
	is a string.
	
	\begin{table}[h!]
	\begin{tabular}{l|ccc|c|}
		& 0 & 1	& 2 & M-N'\\ \hline
		$f_{n,R}$ & $f_{n}(N)$ & $f_{n}(N+1)$ & $f_{n}(N+2)$ & $f_{n}(N+M-N')$ \\
		$f_{n,R}\simeq$ & $g(-M)$ & $g(-M-1)$ & $g(-M-2)$ & $g(-2M+N')$ \\ \hline
		$\overline{g_L}$ & $g(-M)$ & $g(-M-1)$ & $g(-M-2)$ & $g(-2M+N')$ \\
		$\overline{g_L}$ & $g(M)$ & $g(M+1)$ & $g(M+2)$ & $g(2M-N')$  \\ \hline
		$g_R$ & $g(N')$ & $g(N'+1)$ & $g(N'+2)$ & $g(M)$\\ \hline
		$\overline{g_{1,L}}\simeq$ & $g(N')$ & $g(N'+1)$ & $g(N'+2)$ & $g(M)$ \\
		$\overline{g_{1,L}}$ & $g_1(-M')$ & $g_1(-M-1)$ & $g_1(-M'-2)$ & $g_1(-M'-M+N')$ \\ \hline	
	\end{tabular}
	\caption{$\simeq$ means they are semi-coarse homotopic to that row}
	\label{tab:StarOperation}
	\end{table}
		
	\begin{table}[h!]
	\begin{tabular}{l|ccc|}
		& 0 & 1 & 2\\ \hline
		$f_{n,R'}$ & $f_n(N)$ & $f_n(N+1)$ & $f_n(N+2)$ \\
		$f_{n,R'}\simeq$ & $g(M)$ & $g(M+1)$ & $g(M+2)$\\ \hline
		$g_{1,L'}\simeq$ & $g(M)$ & $g(M+1)$ & $g(M+2)$\\
		$g_{1,L'}$ & $g_1(-M-M'+N')$ & $g_1(-M-M'+N'-1)$ & $g_1(-M-M'+N'-2)$\\ \hline
	\end{tabular}
	\caption{$\simeq$ means they are semi-coarse homotopic to that row}
	\label{tab:StarOperationCon}
	\end{table}
\end{observacion}

	Now, with the previous remark, we can define the concatenation between two strings.

\begin{definicion}
\label{def:StarOperationStrings}
	Let $F=(f_1,\ldots,f_n)$ be a path from $\langle f \rangle_\infty$ to $\langle g \rangle_\infty$  and $G=(g_1,\ldots,g_m)$ be a path from $\langle g \rangle_\infty$ to $\langle h \rangle_\infty$.
	We define $(f_1,\ldots,f_n)\star (g_1,\ldots g_m)\coloneqq (f_1,\ldots,f_n,g_1,\ldots,g_m)$.

$\overline{F}\coloneqq (\overline{f_n},\ldots, \overline{f_1})$ is called the opposite direction of $F$.
\end{definicion}

	Although we have defined an operation between pairs of strings, we do not quite have all the elements necessary to obtain a groupoid from this construction.
	At this point, it is not yet clear how to define the identity of a given $\langle f \rangle_\infty$, how to identify $\overline{F}$ with the inverse of $F$, or how the product of an element with its inverse element produces such an identity.
	To achieve these goals, we need a rule that enable us to relate two strings whose difference is a pair of consecutive opposite maps (\autoref{def:DeleteOppositeStrings}).

\begin{definicion}
\label{def:DeleteOppositeStrings}
	Let $n\geq 3$, $i\in \{1,\ldots, n-1\}$ and $F=(f_1,\ldots, f_n)\in \mathcal{F}$.
	Suppose that there exists $N\in\Z$ such that $f_{i}(z) = \overline{f_{i+1}}(z+N)$ for every $z\in \Z$.
	We say that \emph{the result of deleting two consecutive opposite maps in $i$} is $F'=(f_1,\ldots,f_{i-1},f_{i+2},\ldots, f_n)$, i.e. the string $F$ without the maps $f_i$ and $f_{i+1}$.
	We express this relation as $F\xrightarrow{d_{op}(i)} F'$.
	
	In order to simplify future notation, we are going to denote $F\xleftrightarrow{d_{op}(i)} F'$ if $F\xrightarrow{d_{op}(i)} F'$ or $F' \xrightarrow{d_{op}(i)} F$.
	This simplification does not produce ambiguity since one will be a $m$-string, while the other will be a $(m-2)$-string.
	\nomenclature{$\xrightarrow{d_{op}(i)}$}{Delete the opposite maps $f_i$ and $f_{i+1}$ in a string $F$}
\end{definicion}

	We have to verify that the result of deleting two consecutive maps is indeed a string connecting the same two classes of $F$.
	Let $F=(f_1,\ldots,f_n)$ a string from $\langle f_0\rangle_\infty$ to $\langle f_0\rangle_\infty$ and $n$, $i$, $N$, $f_i$ and $f_{i+1}$ as the theorem.
	
	Since $F$ is a string, there exists $N_{i-1},M_{i-1},N_{i+1},M_{i+1}$ such that
	\[f_{i-1,R}(z) \coloneqq f_{i-1}(z+N_{i-1}),\ f_{i,L}(z) \coloneqq f_{i}(z+M_{i-1}) \]
	\[f_{i+1,R}(z) \coloneqq f_{i+1}(z+N_{i+1}),\ f_{i+2,L}(z) \coloneqq f_{i+2}(z+M_{i+1}) \]
	satisfying that $f_{i-1,R}\mid_{[0,\infty)} \simeq_{sc} \overline{f_{i,L}}\mid_{[0,\infty)}$ and $f_{i+1,R}\mid_{[0,\infty)} \simeq_{sc} \overline{f_{i+2,L}}\mid_{[0,\infty)}$.
	Given that $f_{i}(z) \simeq_{sc} \overline{f_{i+1}}(z+N)$, then
	\[\overline{f_{i,L}}(z) = \overline{f_{i}}(z-M_{i-1}) \simeq_{sc} f_{i+1} (z-M_{i-1}-N) = f_{i+1,R} (z-M_{i-1}-N-N_{i+1}). \]
	We have then two cases depending on the integer number $A\coloneqq M_{i-1}+N+N_{i+1}$:
	\begin{itemize}
		\item If $A \leq 0$, then $f_{i-1,R}\mid_{[0,\infty)} \simeq_{sc} \overline{f_{i+2,L}}\mid_{[-A,\infty)}$.
		Thus, it is enough to define $f'_{i+2,L}(z) = f_{i+2,L}(z-A)$ to conclude the desired result.
		\item If $A > 0$, then $f_{i-1,R}\mid_{[A,\infty)} \simeq_{sc} \overline{f_{i+2,L}}\mid_{[0,\infty)}$.
		Thus, it is enough to define $f'_{i-1,R}(z) = f_{i-1,R}(z+A)$ to conclude the desired result
		.
	\end{itemize}
	
	The above \autoref{def:DeleteOppositeStrings} let us add and delete pairs of ``opposite" functions in a string. However, we will be interested in several other operations on strings as well.
	We define one more of them in this section, and let the other for the homotopy in the strings.
	
	Looking for simplify some posterior computations, we also relate a pair of consecutive maps $f$ and $g$ when the right tail of $f$ and the left tail of $g$ are the same constant (\autoref{def:MergeConstants}).

\begin{definicion}
\label{def:MergeConstants}
	Let $F=(f_1, \ldots, f_n)$ a string in $X$ and $x_0\in X$. Suppose that there exists $i\in \{1,\ldots, n-1\}$ and $N,M\in \Z$ such that $f_i(x)=f_{i+1}(y)=x_0$ if $x\geq N$ and $y\leq M$.
	For any $j\in\Z$ such that $N\leq j \leq M$, note that $F'\coloneqq (f_1, \ldots, f_{i-1}, g, f_{i+1}, \ldots, f_n)$
	 where
	\[g(x)\coloneqq \begin{cases}
		f_i(x) & x\leq j, \\
		f_{i+1}(x) & x>j. \end{cases}\]
	is a string. Furthermore, suppose that $g$ satisfies the following two implications (possibly vacuously):
	\begin{itemize}
		\item If $g\mid_{[j-1,\infty)}$ is a periodic map and there exists a positive number $k>1$ such that $g\mid_{[j-k,\infty)}$, then $g\mid_{[j,\infty)}$ is homotopic to a constant.
		
		\item If $g\mid_{(\infty,j+1]}$ is a periodic map and there exists a positive number $k>1$ such that $g\mid_{(\infty,j+k]}$, then $g\mid_{(\infty,j]}$ is homotopic to a constant.
	\end{itemize} 
	
	In this case, we call $F'$ the \emph{the merge of F in the coordinate $i$ at point $j$}, and if a string $F'$ may be obtained from the string $F$ by such a merge, then we write $F\xrightarrow{m(i,j)} (f_1, \ldots, f_{i-1}, g, f_{i+1}, \ldots, f_n)$

	In order to simplify future notation, we are going to denote $F\xleftrightarrow{m(i,j)} F'$ if $F\xrightarrow{m(i,j)} F'$ or $F' \xrightarrow{m(i,j)} F$.
	This simplification does not produce ambiguity since one will be a $m$-string, while the other will be a $(m-1)$-string.
	\nomenclature{$\xrightarrow{m(i,j)}$}{Merging maps $f_i$ and $f_{i+1}$ in a string $F$ which are eventually constant at the element $j$}
\end{definicion}

\begin{observacion}
	The reason for the conditions on $g$ is not yet clear in this section.
	However, these conditions prevent pathological behavior when the equivalence relation $\simeq_d$ (\autoref{def:DelEqRel}) is added.
	
	In simple terms, the conditions state that we can only merge a pair of consecutive functions if we can guarantee that neither of the two tails has an infinite number of cycles.
\end{observacion}

\begin{ejemplo}
\label{Exam:ThreeTurnsInOne}
	We cannot merge every consecutive pair of functions in a string.
	For example, let $C_4$ the semi-coarse space induced for the following graph, called the $4$-cycle graph
	\begin{center}
	\begin{tikzpicture}
		\coordinate (n) at (0,1);
		\coordinate (s) at (0,-1);
		\coordinate (w) at (-1,0);
		\coordinate (e) at (1,0);	

		\filldraw[black] (n) circle (0pt) node[anchor=south]{$3$};
		\filldraw[black] (s) circle (0pt) node[anchor=north]{$1$};
		\filldraw[black] (w) circle (0pt) node[anchor=east]{$0$};
		\filldraw[black] (e) circle (0pt) node[anchor=west]{$2$};

		\draw[line width=0.5mm,black] (n)--(e);
		\draw[line width=0.5mm,black] (n)--(w);

		\draw[line width=0.5mm,black] (s)--(e);
		\draw[line width=0.5mm,black] (s)--(w);

		\filldraw[black] (n) circle (2pt);
		\filldraw[black] (s) circle (2pt);
		\filldraw[black] (w) circle (2pt);
		\filldraw[black] (e) circle (2pt);
	\end{tikzpicture}
\end{center}
	and define $e_1:\Z \to C_4$ as follows (see \autoref{fig:ProductThreeOneTurn})
	\[e_1(z) \coloneqq \begin{cases} z & z\in \{0,1,2,3\}\\ 0 & \text{anywhere else} \end{cases} \]
	We cannot merge $e_1$ and $e_1$ in \autoref{def:MergeConstants} with 
	However, if we define $e'_1$ and $e''_1$ such that
	\[ e'_1(z)=e_1(z-4),\ e''_1(z)=e_1(z-8),\]
	
	we can merge $e_1$ and $e'_1$ and then $e''_1$, obtaining the same result as merging $e'_1$ and $e''_1$ and then $e_1$, in the string $(e_1,e'_1,e''_1)$ (see \autoref{fig:ProductThreeOneTurn}.
	That is
	\[ (e_1,e'_1,e''_1) \xrightarrow{m(1,4)} (e_2,e''_1) \xrightarrow{m(1,8)} (e_3),\ (e_1,e'_1,e''_1) \xrightarrow{m(2,8)} (e_1,e'_2) \xrightarrow{m(1,4)} (e_3) \]
	with the maps $e_2$, $e'_2$ and $e_3$ defined as follows (see \autoref{fig:e2cyclegraph} and \autoref{fig:e3cyclegraph}, respectively)
	\[e_2(z) \coloneqq \begin{cases} 1 & z\in \{1,5\},\\ 2 & z\in \{2,6\},\\ 3 & z \in \{3,7\},\\ 0 & \text{anywhere else}, \end{cases}\ e'_2(z) \coloneqq \begin{cases} 1 & z\in \{5,9\},\\ 2 & z\in \{6,10\},\\ 3 & z \in \{7,11\},\\ 0 & \text{anywhere else}, \end{cases}\ e_3(z) \coloneqq \begin{cases} 1 & z\in \{1,5,9\},\\ 2 & z\in \{2,6,10\},\\ 3 & z \in \{3,7,11\},\\ 0 & \text{anywhere else}. \end{cases} \]
	Observe that we mix a string of three maps with one turn in a string of just one map with three consecutive turns.
	
	It would be desirable to merge $e_1$ and $e_1$. We resolve this issue with homotopy when we define the fundamental groupoid (\autoref{def:FundamentalGroupoid}).

\end{ejemplo}

\begin{figure}[h]
	\begin{tikzpicture}
	\foreach \y in {0,1,2}
		{\foreach \x in {-3,-2,-1}
			{\filldraw[black] (\x/6+10/3*\y,0) circle (0.2pt);}
		\foreach \x in {0,1,2,6,7,8}
			{
			\filldraw[RedOrange] (\x/3+10/3*\y,0) circle (1.5pt) node[below,black] {\tiny{\inteval{\x-2+4*\y}}};}
		\filldraw[RedOrange] (1+10/3*\y,1/3) circle (1.5pt);
		\filldraw[Gray!30] (1+10/3*\y,0) circle (1.5pt) node[below,black] {\tiny{\inteval{1+4*\y}}};
		\filldraw[RedOrange] (4/3+10/3*\y,2/3) circle (1.5pt);
		\filldraw[Gray!30] (4/3+10/3*\y,0) circle (1.5pt) node[below,black] {\tiny{\inteval{2+4*\y}}};
		\filldraw[RedOrange] (5/3+10/3*\y,1) circle (1.5pt);
		\filldraw[Gray!30] (5/3+10/3*\y,0) circle (1.5pt) node[below,black] {\tiny{\inteval{3+4*\y}}};}
		\foreach \x in {-3,-2,-1}
			{\filldraw[black] (\x/6+30/3,0) circle (0.2pt);}
\end{tikzpicture}
\caption{Representation of $(e_1,e'_1,e''_1)$}
\label{fig:ProductThreeOneTurn}
\end{figure}

\begin{figure}[h]
	\begin{tikzpicture}
		\foreach \x in {-3,-2,-1}
			{\filldraw[black] (\x/6,0) circle (0.2pt);}
		\filldraw[RedOrange] (0,0) circle (1.5pt) node[below,black] {\tiny{-1}};
		\foreach \x in {0,1}
			{\filldraw[RedOrange] (1/3+4/3*\x,0) circle (1.5pt) node[below,black] {\tiny{\inteval{4*\x}}};
			\filldraw[RedOrange] (2/3+\x*4/3,1/3) circle (1.5pt);
			\filldraw[Gray!30] (2/3+\x*4/3,0) circle (1.5pt) node[below,black] {\tiny{\inteval{4*\x+1}}};
			\filldraw[RedOrange] (1+\x*4/3,2/3) circle (1.5pt);
			\filldraw[Gray!30] (1+\x*4/3,0) circle (1.5pt) node[below,black] {\tiny{\inteval{4*\x+2}}};
			\filldraw[RedOrange] (4/3+\x*4/3,1) circle (1.5pt);
			\filldraw[Gray!30] (4/3+\x*4/3,0) circle (1.5pt) node[below,black] {\tiny{\inteval{4*\x+3}}};}
		\filldraw[RedOrange] (3,0) circle (1.5pt) node[below,black] {\tiny{8}};
		\foreach \x in {-3,-2,-1}
			{\filldraw[black] (\x/6+11/3,0) circle (0.2pt);}
\end{tikzpicture}\hspace{5mm}
	\begin{tikzpicture}
		\foreach \x in {-3,-2,-1}
			{\filldraw[black] (\x/6,0) circle (0.2pt);}
		\filldraw[RedOrange] (0,0) circle (1.5pt) node[below,black] {\tiny{3}};
		\foreach \x in {0,1}
			{\filldraw[RedOrange] (1/3+4/3*\x,0) circle (1.5pt) node[below,black] {\tiny{\inteval{4*\x+4}}};
			\filldraw[RedOrange] (2/3+\x*4/3,1/3) circle (1.5pt);
			\filldraw[Gray!30] (2/3+\x*4/3,0) circle (1.5pt) node[below,black] {\tiny{\inteval{4*\x+5}}};
			\filldraw[RedOrange] (1+\x*4/3,2/3) circle (1.5pt);
			\filldraw[Gray!30] (1+\x*4/3,0) circle (1.5pt) node[below,black] {\tiny{\inteval{4*\x+6}}};
			\filldraw[RedOrange] (4/3+\x*4/3,1) circle (1.5pt);
			\filldraw[Gray!30] (4/3+\x*4/3,0) circle (1.5pt) node[below,black] {\tiny{\inteval{4*\x+7}}};}
		\filldraw[RedOrange] (3,0) circle (1.5pt) node[below,black] {\tiny{12}};
		\foreach \x in {-3,-2,-1}
			{\filldraw[black] (\x/6+11/3,0) circle (0.2pt);}
\end{tikzpicture}
\caption{Representation of $e_2$ and $e'_2$, respectively}
\label{fig:e2cyclegraph}
\end{figure}

\begin{figure}[h]
	\begin{tikzpicture}
		\foreach \x in {-3,-2,-1}
			{\filldraw[black] (\x/6,0) circle (0.2pt);}
		\filldraw[RedOrange] (0,0) circle (1.5pt) node[below,black] {\tiny{-1}};
		\foreach \x in {0,1,2}
			{\filldraw[RedOrange] (1/3+4/3*\x,0) circle (1.5pt) node[below,black] {\tiny{\inteval{4*\x}}};
			\filldraw[RedOrange] (2/3+\x*4/3,1/3) circle (1.5pt);
			\filldraw[Gray!30] (2/3+\x*4/3,0) circle (1.5pt) node[below,black] {\tiny{\inteval{4*\x+1}}};
			\filldraw[RedOrange] (1+\x*4/3,2/3) circle (1.5pt);
			\filldraw[Gray!30] (1+\x*4/3,0) circle (1.5pt) node[below,black] {\tiny{\inteval{4*\x+2}}};
			\filldraw[RedOrange] (4/3+\x*4/3,1) circle (1.5pt);
			\filldraw[Gray!30] (4/3+\x*4/3,0) circle (1.5pt) node[below,black] {\tiny{\inteval{4*\x+3}}};}
		\filldraw[RedOrange] (13/3,0) circle (1.5pt) node[below,black] {\tiny{12}};
		\foreach \x in {-3,-2,-1}
			{\filldraw[black] (\x/6+15/3,0) circle (0.2pt);}
\end{tikzpicture}
\caption{Representation of $e_3$}
\label{fig:e3cyclegraph}
\end{figure}

	We form an equivalence relation with these two ingredients.
	The specific relation and result are established in the following sentences:
 
\begin{definicion}
	Let $F,G\in \mathcal{F}$.
	We say that \emph{we can convert $F$ in $G$} if there exists a finite sequence of arrows
\begin{align*}
F\coloneqq F_0 \xleftrightarrow{\alpha_0} F_1 \xleftrightarrow{\alpha_1} \cdots \xleftrightarrow{\alpha_{n-1}} F_n \eqqcolon G
\end{align*}
with $\alpha_k\in \{ d_{op}(i) \mid i\in\mathbb{Z} \}\cup \{ m(i,j) \mid i,j\in\mathbb{Z} \}$.
	We express this relation as $F\simeq_S G$.
	\nomenclature{$\simeq_S$}{The equivalence relation formed by the relations merge and delete opposite maps}
\end{definicion}

\begin{teorema}
\label{theo:SisEqRel}
	$\simeq_S$ is an equivalence relation.
\end{teorema}

\begin{proof}
	Reflexivity and symmetry are directly derived from the definition because
	\[ (f) \xleftarrow{d_{op}(i)} (f,\overline{f},f) \xrightarrow{d_{op}(i)} (f) ,\]
	and the relations $\xleftrightarrow{d_{op}(i)}$ and $\xleftrightarrow{m(i)}$ are already symmetric, respectively.
	
	For transitivity, let $F,G,H\in \mathcal{F}$ with $F=(f_1,\ldots, f_n)$.
	Suppose that $F\simeq_S G$ and $G\simeq_S H$. Then, there exist the sequence
\begin{align*}
F\coloneqq F_0 \xleftrightarrow{\alpha_0} F_1 \xleftrightarrow{\alpha_1} \cdots \xleftrightarrow{\alpha_{n-1}} F_n \eqqcolon G\\
G\coloneqq G_0 \xleftrightarrow{\beta_0} G_1 \xleftrightarrow{\beta_1} \cdots \xleftrightarrow{\beta_{n-1}} G_n \eqqcolon H
\end{align*}
with $\alpha_i\in \{ d_{op}(i) \mid i\in\mathbb{Z} \}\cup \{ m(i) \mid i\in\mathbb{Z} \}\cup \{ = \}$. Therefore
\begin{align*}
F = F_0 \xleftrightarrow{\alpha_0} F_1 \xleftrightarrow{\alpha_1} \cdots \xleftrightarrow{\alpha_{n-1}} F_n = G = G_0 \xleftrightarrow{\beta_0} G_1 \xleftrightarrow{\beta_1} \cdots \xleftrightarrow{\beta_{n-1}} G_n = H
\end{align*}
and thus $F\simeq_S H$.
\end{proof}

	Let $F=(f_1,\ldots,f_n) \in \mathcal{F}$, we denote the class of $F$ as $\langle F\rangle_S$ or $\langle f_1,\ldots, f_n\rangle_S$.
	\nomenclature{$\langle F\rangle_S$}{The equivalence class of $F$ with relation $\simeq_S$} 		Finally, we define the concatenation between classes and see that is well-defined.

\begin{proposicion}
\label{prop:ProductClassSWellDef}
Let $F,G\in\mathcal{F}$. Then $\langle F \rangle_S\star \langle G \rangle_S = \langle F\star G \rangle_S$ is well-defined.
\end{proposicion}

\begin{proof}
Let $F,G\in\mathcal{F}$, $F'\in \langle F \rangle_S$ and $G'\in \langle G \rangle_S$, then there exist the sequences
\begin{align*}
F \coloneqq F_0 \xrightarrow{\alpha_0} F_1 \xrightarrow{\alpha_1} \cdots \xrightarrow{\alpha_{n-1}} F_n \eqqcolon F',\\
G \coloneqq G_0 \xrightarrow{\beta_0} G_1 \xrightarrow{\beta_1} \cdots \xrightarrow{\beta_{n-1}}  G_m \eqqcolon G'.
\end{align*}
Then we have that
\begin{align*}
F \star G = & F_0\star G \xleftrightarrow{\alpha_0\star 1_G} F_1\star G \xleftrightarrow{\alpha_1\star 1_G} \cdots \xleftrightarrow{\alpha_{n-1}\star 1_G} F_n \star G= F'\star G,\\
F' \star G = & F'\star G_0 \xleftrightarrow{1_{F'}\star \beta_{0}} F'\star G_1 \xleftrightarrow{1_{F'}\star \beta_{1}} \cdots \xleftrightarrow{1_{F'}\star \beta_{m-1}} F'\star G_m = F'\star G'.\qedhere
\end{align*}
.
\end{proof}

	Considering this class, we observe in the following lemma that with this operation every element has  inverses and left and right identities.
	 Hence we are able to define the groupoid. 

\begin{lema}
\label{lem:InverseExtendedGroupoid}
Let $F=(f_1, \ldots, f_n)\in\mathcal{F}$, then $\langle F \rangle_S \star \langle \overline{F} \rangle_S= \langle f_1,\overline{f_1} \rangle_S$ and $\langle \overline{F} \rangle_S \star \langle F \rangle_S = \langle \overline{f_n},f_n \rangle_S$. Moreover, $\langle f_1,\overline{f_1} \rangle_S \star \langle F \rangle_S = \langle F \rangle_S = \langle F \rangle_S \star \langle \overline{f_n},f_n \rangle_S$.
\end{lema}

\begin{proof}
Let $F=(f_1, \ldots, f_n)\in\mathcal{F}$. Suppose that $n>1$, then
\begin{align*}
	\langle F \rangle_S \star \langle \overline{F} \rangle_S = & \langle f_1,\ldots, f_n \rangle_S \star \langle \overline{f_n},\ldots, \overline{f_1} \rangle_S = \langle f_1,\ldots, f_n,\overline{f_n},\ldots, \overline{f_1} \rangle_S\\
 = & \langle f_1,\dots,f_{n-1},\overline{f_{n-1}},\ldots,\overline{f_1} \rangle_S = \langle f_1,\overline{f_1} \rangle_S,\\ 
	\langle \overline{F} \rangle_S \star \langle F \rangle_S = & \langle \overline{f_n},\ldots, \overline{f_1} \rangle_S \star \langle f_1,\ldots, f_n \rangle_S = \langle \overline{f_n},\ldots, \overline{f_1},f_1,\ldots, f_n \rangle_S\\
 = & \langle \overline{f_n},\ldots, \overline{f_2},f_2,\ldots, f_n \rangle_S = \langle \overline{f_n},f_n \rangle_S,\\
	\langle f_1,\overline{f_1} \rangle_S \star \langle F \rangle_S = & \langle f_1,\overline{f_1} \rangle_S \star \langle f_1,\ldots,f_n \rangle_S = \langle f_1,\overline{f_1},f_1,\ldots f_n \rangle_S\\
 = & \langle f_1,\ldots,f_n \rangle_S = \langle F \rangle_S\\
	\langle F \rangle_S \star \langle \overline{f_n},f_n \rangle_S = & \langle f_1,\ldots,f_n \rangle_S \star \langle \overline{f_n},f_n \rangle_S = \langle f_1,\ldots,f_n,\overline{f_n},f_n \rangle_S\\
  = & \langle f_1,\ldots,f_n \rangle_S = \langle F \rangle_S \qedhere
\end{align*}
\end{proof}

\begin{lema}
\label{lem:SameIdentity}
	Let $f,g:\mathbb{Z}_1\rightarrow X$ be bornologous maps such that $( \overline{f}, g )$ is a string. Then, $\langle f, \overline{f} \rangle_S = \langle g, \overline{g} \rangle_S$.
	In particular, if $f$ is a symmetric map and $g \in \langle f \rangle_\infty$,
	then $\langle f, \overline{f} \rangle_S = \langle g, \overline{g} \rangle_S$ and $\langle \overline{f}, f \rangle_S = \langle \overline{g}, g \rangle_S$,
	i.e. the object $\langle f \rangle_\infty$ have left and right identities.
\end{lema}

\begin{proof}
	Given $f,g$ such that $(\overline{f},g)$ is a string, we can operate $(f)\star (\overline{f},g) \star (\overline{g})$ because both pair $(f,\overline{f})$ and $(g,\overline{g})$ are opposite.
	Then, by definition, we have that
	\[( f, \overline{f}) \xleftarrow{d_{op}(3)} ( f, \overline{f}, g, \overline{g}) \xrightarrow{d_{op}(1)} ( g, \overline{g} ).\]
	Thus $\langle f, \overline{f} \rangle_S = \langle g, \overline{g} \rangle_S$,
	obtaining the first part of the desired result.
	
	In brief, observe that $(\overline{f},g)$ and $(f, \overline{g})$ are strings if $f$ is symmetric and $g\in \langle g\rangle_\infty$.
	Concluding the second part of the lemma.
\end{proof}

\begin{proposicion}
\label{prop:Identity1String}
	Let $F=(f_1,\ldots,f_n)$ be a $n$-string from $\langle f \rangle_\infty$ to $\langle g \rangle_\infty$.
	Then $(f)\star F = F = F\star (g)$.
	Thus, we can conclude that $(f')\star F = F = F\star (g')$ if $f'\in\langle f \rangle_\infty$ and $g'\in\langle g\rangle_\infty$.
\end{proposicion}

\begin{proof}
	Let $F,f,g$ be as the statement.
	Define $f^{\leftarrow}, f^{\rightarrow}, g^{\leftarrow}, g^{\rightarrow}:\Z\to X$ as follows:
	\[f^{\leftarrow}(z) \coloneqq \begin{cases}
	f(z) & z\leq 0,\\ f(0) & z>0,
	\end{cases}\hspace{5mm} f^{\rightarrow}(z) \coloneqq \begin{cases}
	f(z) & z\geq 0,\\ f(0) & z<0,
	\end{cases}\]
	\[g^{\leftarrow}(z)\coloneqq \begin{cases}
	g(z) & z\leq 0,\\ g(0) & z>0,
	\end{cases}\hspace{5mm} g^{\rightarrow}(z)\coloneqq \begin{cases}
	g(z) & z\geq 0,\\ g(0) & z<0.
	\end{cases} \]
	Observe that $f^\leftarrow = \overline{f^\rightarrow}$ and $g^\leftarrow = \overline{g^\rightarrow}$.
	In addition, we have that
	\[ (f^\leftarrow,f^\rightarrow) \xrightarrow{m(1,0)} (f),\ (g^\leftarrow,g^\rightarrow) \xrightarrow{m(1,0)} (g).\]
	By \autoref{lem:SameIdentity}, we conclude that
	\[ \langle (f)\star F \rangle_S = \langle (f^\leftarrow, f^\rightarrow)\star F \rangle_S = \langle F \rangle_S = \langle F \star (g^\leftarrow, g^\rightarrow) \rangle_S = \langle F \star (g) \rangle_S .\qedhere\]
\end{proof}	
	
	Given the \autoref{def:StarOperationStrings}, we automatically have that the product of two strings is associative. Now, with \autoref{lem:InverseExtendedGroupoid} and \autoref{lem:SameIdentity}, we have that every string has a left inverse and a right inverse; in addition, every string has a right identity and a left identity.
	In addition, \autoref{prop:Identity1String} claims that $\langle f\rangle_\infty$ is the left identity and $\langle g \rangle_\infty$ is the right identity of $F$ if $F$ is a string from $\langle f\rangle_\infty$ to $\langle g \rangle_\infty$.
	We assign a name to this groupoid.

\begin{definicion}
\label{def:ExtendGrupoid}
The \emph{extended groupoid of $X$}, $\overline{\pi}_{\leq 1}(X)$, is the groupoid with objects the symmetric maps to $X$ and morphisms the strings.\nomenclature{$\overline{\pi}_{\leq 1}(X)$}{Extended groupoid of $X$}
\end{definicion}

	Suppose that we have $X$ and $Y$ semi-coarse spaces and a bornologous map $h:X\rightarrow Y$.
	On the one hand, if we have a string $(f_1,\ldots,f_n)$ on $X$, then $h(f_1,\ldots,f_n)\coloneqq (hf_1,\ldots,hf_n)$ is a string on $Y$.
	On the other hand, if we have a symmetric map $f$, then $hf$ is still symmetric, and $hf\simeq_{sc} hg$ if $f\simeq_{sc} g$.
	That is, in this section, we have developed a functor from the category of semi-coarse spaces to the category of groupoids.
	In particular, we directly observe that the extended groupoid is a semi-coarse invariant.
	In addition, if there exists a bornologous map $f:X\to Y$ and $g:Y\to X$ such that $gf\simeq_{sc}\Id_X$ and $fg \simeq_{sc} \Id_Y$, then we obtain an isomorphism between the objects of both categories.

\begin{proposicion}
\label{prop:ExtGroupoidSCInv}
	Let $X,Y$ be isomorphic semi-coarse spaces.
	Then $\overline{\pi}_{\leq 1}(X)$ and $\overline{\pi}_{\leq 1}(Y)$ are isomorphic as groupoids.
	In particular, if $X,Y$ are strong homotopy equivalent, then there exists an isomorphism between the objects of $\overline{\pi}_{\leq 1}(X)$ and the objects of $\overline{\pi}_{\leq 1}(Y)$
\end{proposicion}

\begin{corolario}
\label{cor:ExtGroupoidSCInv}
	Any semi-coarse space $X$ such that the objects of $\overline{\pi}_{\leq 1}(X)$ has cardinality greater than $1$ is not strong homotopy equivalent to $X_\text{ind}$.
\end{corolario}

\begin{corolario}	
	Any semi-coarse space $X$ with $n$ connected components has at least $n$ objects in $\overline{\pi}_{\leq 1}(X)$.
\end{corolario}

\section{Semi-coarse Fundamental Groupoid}
\label{cap.FundGrupoid}

	In the previous \autoref{cap.SCStrings}, we define a small category associated to a every semi-coarse space $X$ giving a class of equivalence in the strings of $X$.
	Although this small category is actually a groupoid, it is not easy to observe what kind of information is captured for this groupoid.
	The massive amount of arrows in this category only allows to observe that every arrow has an inverse and that the semi-coarse space is somehow embedded in the extended groupoid of $X$.
	
	In this section we increase the number of relations between strings by one compared with the extended groupoid.
	This relation translates the continuous homotopy in topological spaces to the discrete context of semi-coarse spaces.
	The relation and its philosophy is the following:
	\begin{itemize}		
		\item Delete-one-point equivalence, in \autoref{def:DeletingAndAdding}, mimics the topological  homeomorphisms from $[0,1]$ to $[0,2]$, allowing us to ``create'' a finite number of values between two consecutive integers.
	\end{itemize}

	We develop this notion and conclude the chapter by defining the semi-coarse fundamental groupoid considering this relation plus those in \autoref{def:DeleteOppositeStrings} and \autoref{def:MergeConstants}.

	In comparison with the semi-coarse homotopy, the following relation only relates two maps which are homotopy equivalent in a bounded discrete interval of the form $[-N,N]$ (\autoref{prop:FiniteHomotopyIFFDeleteEq}).
	However, this new relation allows to add a point ``between'' two consecutive integers, which is in general technical impossible just using the semi-coarse homotopy and essential for the proof of \autoref{lem:LebesgueCoveringSC}.
	
\begin{definicion}
\label{def:DeletingAndAdding}
	Let $f:\mathbb{Z}_1\rightarrow X$ be a bornologous map.
	Let $z_0\in\mathbb{Z}$ and $x_0\in X$.
	\begin{itemize}
		\item We said that $g$ is the result of deleting $z_0$ in $f$, or we can delete $z_0$ in $f$, if 
\begin{align*}
g(z)\coloneqq \left\lbrace\begin{array}{ll}
f(z) & \text{ if }z<z_0\\
f(z+1) & \text{ if }z\geq z_0.
\end{array}\right.
\end{align*}
		is bornologous and $g$ satisfies that
		\begin{itemize}
			\item If $g\mid_{[j-1,\infty)}$ is a periodic map and there exists a positive number $k>1$ such that $g\mid_{[j-k,\infty)}$, then $g\mid_{[j,\infty)}$ is homotopic to a constant.
		
			\item If $g\mid_{(\infty,j+1]}$ is a periodic map and there exists a positive number $k>1$ such that $g\mid_{(\infty,j+k]}$, then $g\mid_{(\infty,j]}$ is homotopic to a constant.
		\end{itemize}
		This is denoted as $f\xrightarrow{d(z_0)} g$.
		\nomenclature{$\xrightarrow{d(z_0)}$}{In a bornologous map $f:\Z \to X$, delete the coordinate $z_0$ and move all the right elements one step to the left}

		\item We said that $g$ is the result of adding $x_0$ at $z_0$ in $f$, or we can add $x_0$ at $z_0$, if
\begin{align*}
g(z)\coloneqq \left\lbrace\begin{array}{ll}
f(z) & \text{ if }z<z_0\\
x_0 & \text{ if }z=z_0\\
f(z-1) & \text{ if }z> z_0.
\end{array}\right.
\end{align*}
is bornologous and $f$ satisfies that
	\begin{itemize}
		\item If $f\mid_{[j-1,\infty)}$ is a periodic map and there exists a positive number $k>1$ such that $f\mid_{[j-k,\infty)}$, then $f\mid_{[j,\infty)}$ is homotopic to a constant.
		
		\item If $f\mid_{(\infty,j+1]}$ is a periodic map and there exists a positive number $k>1$ such that $f\mid_{(\infty,j+k]}$, then $f\mid_{(\infty,j]}$ is homotopic to a constant.
	\end{itemize}.
		This is denoted as $f\xrightarrow{a(z_0,x_0)} g$.
	\nomenclature{$f\xrightarrow{a(z_0,x_0)} g$}{In a bornologous map $f:\Z \to X$, move to the right all the elements after the coordinate $z_0$ and in the coordinate $z_0$ put $x_0$}
\end{itemize}
\end{definicion}

\begin{lema}
\label{lem:AddAndDeletePoints}
	Let $X$ be a semi-coarse space and $f,g:\mathbb{Z}_1\rightarrow X$ be bornologous maps.
	Then, $f\xrightarrow{d(z_0)} g$ if and only if $g\xrightarrow{a(z_0,f(z_0))} f$.
\end{lema}

\begin{proof}
	Suppose that $f\xrightarrow{d(z_0)} g$.
	By definition, we have that $g(z)=f(z)$ for every $z<z_0$ and $g(z)=f(z+1)$ for every $z\geq z_0$.
	Moreover, the result of adding $f(z_0)$ at $z_0$ in $g$ is
\begin{align*}
h(z) = \left\lbrace\begin{array}{ll}
g(z)   & \text{ if } z<z_0\\
f(z_0) & \text{ if } z=z_0\\
g(z-1)   & \text{ if } z>z_0.
\end{array}\right.
\end{align*}
	Replacing the values of $g$, we have that
\begin{align*}
h(z) = \left\lbrace\begin{array}{ll}
f(z)            & \text{ if } z<z_0\\
f(z_0)          & \text{ if } z=z_0\\
f(z-1+1)=f(z)   & \text{ if } z>z_0.
\end{array}\right.
\end{align*}
	Thus $g\xrightarrow{a(z_0,f(z_0))} f$.

	Suppose that $g\xrightarrow{a(z_0,f(z_0))} f$.
	By definition we have that $f(z)=g(z)$ for every $z<z_0$ and $f(z)=g(z-1)$ for every $z>z_0$.
	Moreover, the result of deleting $z_0$ in $f$ is
\begin{align*}
h(z) = \left\lbrace\begin{array}{ll}
f(z)     & \text{ if } z<z_0\\
f(z+1)   & \text{ if } z\geq z_0.
\end{array}\right.
\end{align*}
	Replacing the values of $f$, we have that
\begin{align*}
h(z) = \left\lbrace\begin{array}{ll}
g(z)            & \text{ if } z<z_0\\
g(z+1-1)=g(z)   & \text{ if } z\geq z_0.
\end{array}\right.
\end{align*}
	Thus $f\xrightarrow{a(z_0)} f$.
\end{proof}

\begin{lema}
\label{lem:TrivialAdding}
	Let $X$ be a semi-coarse space and $f:\mathbb{Z}_1\rightarrow X$ be a bornologous map, then
\begin{align*}
f_{a(z_0)}\coloneqq \left\lbrace\begin{array}{ll}
f(z)   & \text{ if } z<z_0\\
f(z_0) & \text{ if } z=z_0\\
f(z-1)   & \text{ if } z>z_0
\end{array}\right.
\end{align*}
	is a bornologous map, and therefore $f \xrightarrow{a(z_0,f(z_0))}f_{a(z_0)}$ whenever $f$ and $z_0$ satisfy the conditions of \autoref{def:DeletingAndAdding}.
\end{lema}

\begin{proof}
	Let $f:\mathbb{Z}_1\rightarrow X$ be a bornologous maps.
	We take the subsets $X_1\coloneqq (-\infty, z_0]$, $X_2\coloneqq [z_0,z_0+1]$ and $X_3 \coloneqq [z_0+1,\infty)$.
	Observe that every controlled set in $\mathbb{Z}_1$ is an union of controlled sets in $\mathbb{Z}_1\mid_{X_1}$, $\mathbb{Z}_1\mid_{X_2}$ and $\mathbb{Z}_1\mid_{X_3}$, and $X= X_1\cup X_2\cup X_3$.
	Thus, we can use \autoref{prop:BornologousPartition} to prove that $f_{a(z_0)}$ is a bornologous map.
	Observe that $f_{a(z_0)}\mid_{X_1} = f\mid_{X_1}$ and $f_{a(z_0)}\mid_{X_3} = f\mid_{[z_0,\infty)}$, which are bornologous maps because they are restrictions of a bornologous map.
	Further, $f_{a(z_0)}\mid_{X_2}=f\mid_{\{z_0\}}$ which goes from just one point, then it is a bornologous map.
	Thus, $f_{a(z_0)}$ is a bornologous map.
\end{proof}

\begin{definicion}
\label{def:DelEqRel}
	Let $X$ be a semi-coarse space and $f,g:\mathbb{Z}_1\rightarrow X$ be bornologous maps.
	We say that $f\simeq_d g$ if there exist a finite sequence such that
	\nomenclature{$f\simeq_d g$}{The equivalence relation formed by $\xrightarrow{d(z_0)}$ and $\xrightarrow{a(z_0,x_0)}$}
\begin{align*}
f = f_0 \xrightarrow{\alpha_0} f_1 \xrightarrow{\alpha_1}\cdots \xrightarrow{\alpha_{n-1}} f_n = g
\end{align*}
	with $\alpha_i \in \{ d(z) \mid z\in\mathbb{Z} \}\cup \{ a(z,x) \mid z\in\mathbb{Z},x\in X \}$.
\end{definicion}

\begin{teorema}
\label{theo:DeleteIsEqRelation}
	$\simeq_d$ is an equivalence relation.
\end{teorema}

\begin{definicion}
\label{def:SimdStrings}
	Let $F=(f_1,\ldots, f_n)$ and $G=(g_1,\ldots,g_n)$ be $n$ strings.
	We say that $F \simeq_d G$ if $f_i\simeq_d g_i$.
	\nomenclature{$F\simeq_d G$}{The equivalence relation of string formed by $\xrightarrow{d(z_0)}$ and $\xrightarrow{a(z_0,x_0)}$ in every component of $F$ and $G$}
\end{definicion}

	In the result below we prove that \autoref{def:DelEqRel} allows to make homotopy deformations in the path which joins two tails.
	This proposition is the main argument for what the fundamental groupoid (\autoref{def:FundamentalGroupoid}) contains the semi-coarse fundamental group for every point in the space.

\begin{proposicion}
\label{prop:FiniteHomotopyIFFDeleteEq}
	Let $X$ be a semi-coarse space and $f,g:\mathbb{Z}_1\rightarrow X$ such that there exists $N\in\mathbb{N}$ satisfying
\begin{itemize}
	\item $f\mid_{(-\infty,-N]\cup [N,\infty)}=g\mid_{(-\infty,-N]\cup [N,\infty)}$.
	
	\item $f\mid_{[-N,N]}\simeq_{sc} g\mid_{[-N,N]}$ through $H:[-N,N]\times \mathbb{Z}_1\rightarrow X$ such that $H(-N,-)=f(-N)$ and $H(N,-)=f(N)$.
\end{itemize}
Then $f \simeq_d g $.
\end{proposicion}

\begin{proof}
	Let $X$ be a semi-coarse space and $f,g:\mathbb{Z}_1\rightarrow X$ and $\{0,1\}$ with the subspace structure.
	First we consider a homotopy with just one step.
	Suppose that there exists $N\in\mathbb{N}$ satisfying
\begin{itemize}
	\item $f\mid_{(-\infty,-N]\cup [N,\infty)}=g\mid_{(-\infty,-N]\cup [N,\infty)}$.
	
	\item There exists $H:[-N,N]\times \{0,1\}\rightarrow X$ bornologous map such that $H(-N,-)=f(-N)$, $H(N,-)=f(N)$, $H(z,0)=f(z)$ and $H(z,1)=g(z)$.
\end{itemize}

	By the condition two, observe that $\{ (f(i),g(i+1)), (g(i+1),f(i+2)) \}$ is controlled in $X$ for every $i\in\mathbb{Z}$, in particular for $-N\leq i \leq N$.
	Thus, we obtain the following
\begin{align*}
f \xrightarrow{a(-N+1,g(-N+1))} f'_1 \xrightarrow{d(-N+2)} f''_1.
\end{align*}
	Observe that $f(z)=f''_1(z)$ for every $z\neq -N+1$ and $f''_1(-N+1) = g(-N+1)$.
	We can repeat this interaction as
\begin{align*}
f''_i  \xrightarrow{a(-N+i+1,g(-N+i+1))} f'_{i+1} \xrightarrow{d(-N+i+2)} f''_{i+1},
\end{align*}
	obtaining that $f_i''(z)=f''_{i+1}(z)$ for every $z\neq -N+i+1$ and $f''_{i+1}(-N+i+1)=g(-N+i+1)$.
	Making this algorithm until $i=2N-2$ we obtain $f''_{2N-1}=g$. Thus $f \simeq_d g$.

	In the general case, we have that there exists $N\in\mathbb{N}$ satisfying
\begin{itemize}
	\item $f\mid_{(-\infty,-N]\cup [N,\infty)}=g\mid_{(-\infty,-N]\cup [N,\infty)}$.
	
	\item $f\mid_{[-N,N]}\simeq_{sc} g\mid_{[-N,N]}$ through $H:[-N,N]\times \mathbb{Z}_1\rightarrow X$ such that $H(-N,-)=f(-N)$ and $H(N,-)=f(N)$.
\end{itemize}
	Since $H$ is a homotopy, there exists $a<b\in\mathbb{Z}$ such that $H(z,x)=f(x)$ for every $x\leq a$ and $H(z,x)=g(z)$ for every $x\geq b$.
	Define $f_i:[-N,N]\rightarrow X$ such that $f_i(z)=H(z,i)$ for every $a\leq i \leq b$.
	Hence, taking $H_i \coloneqq H\mid_{[-N,N]\times \{i,i+1\}}$ we obtain that
\begin{align*}
f = f_{a} \simeq_d f_{a+1} \simeq_d \cdots \simeq_{d} f_{b-1} \simeq_{d} f_{b}= g
\end{align*}
	Concluding that $f\simeq_d g$.
\end{proof}

\begin{definicion}
\label{def:SimStringWithMerge}
Let $F,G\in \mathcal{F}$. We say that \emph{we can convert and deform $F$ in $G$} if there exists a finite sequence of arrows
\begin{align*}
F\coloneqq F_0 \simeq_{\alpha_0} F_1 \simeq_{\alpha_1} \cdots \simeq_{\alpha_{n-1}} F_n \eqqcolon G
\end{align*}
	with $\alpha_i\in \{d,S,sc\}$. We express this relation as $F\simeq G$.
	We denote the equivalence class of $F$ by $[F]$.
	\nomenclature{$F\simeq G$}{The equivalence relation of strings formed by $\simeq_S$ and $\simeq_d$}
\end{definicion}

\begin{definicion}
\label{def:FundamentalGroupoid}
	The \emph{fundamental groupoid of $X$}, $\pi_{\leq 1}(X)$, is the groupoid with objects the symmetric maps to $X$ and morphisms the equivalence classes of strings under the relation $\simeq$ defined in \autoref{def:SimStringWithMerge} above.
\end{definicion}

\begin{observacion}
	The fundamental groupoid is actually a groupoid given the relation $\simeq_S$ and the definition of  the $\star$-operation of two strings.
	As we mention in the previous section, $\simeq_S$ provides us of inverses and right and left identities;
	the $\star$-operation ensure that the composition of arrows is defined.
\end{observacion}

	The three relationships defined above give rise to the following two results.
	The first result (\autoref{prop:Affinity}) tells us that any coordinate of a string $F=(f_1,\ldots,f_n)$ can be shifted to the left or right, and the result will still be equivalent to $F$.
	On the other hand, the second result (\autoref{prop:CutEqualTails}) ensures that if the left tail of f is eventually equal to the right tail of g, we can cut both tails and paste the parts of the images which remain.
	
\begin{proposicion}
\label{prop:Affinity}
	Let $X$ be a semi-coarse space and $F=(f_1,\ldots,f_n)$ be a string.
	If we denote by $f'_i:\Z_1\to X$ the map $f'_i(z)=f_i(z-1)$, then
	\[ [F]=[f_1,\ldots,f'_i,\ldots,f_n] .\]
\end{proposicion}

\begin{proof}
	Let $f_i:\Z_1\to X$ be a bornologous maps. There exists a $z_0\in\Z$ such that we can define the maps
	\[ f_i^\leftarrow(z) \coloneqq \begin{cases}
	f_i(z) & z\leq z_0,\\
	f_i(0) & z>z_0,
	\end{cases}\hspace{5mm} f_i^\rightarrow(z) \coloneqq \begin{cases}
	f_i(z) & z\geq z_0,\\
	f_i(0) & z<z_0,
	\end{cases}\]
	and
	\[ (f_i) \xleftarrow{m(1,z_0)} (f_i^\leftarrow,f_i^\rightarrow). \]
	Now define the following maps
	\[ f'_i(z) \coloneqq f_i(z-1),\ g_i^\leftarrow(z) \coloneqq f_i^\leftarrow(z-1),\ g_{i}^\rightarrow(z) \coloneqq f_{i+1}^\rightarrow(z-1).\]
	Thus we have the following relations
	\[ (f_i^\leftarrow,f_i^\rightarrow) \xleftarrow{d_{op}(2)} (f_i^\leftarrow,\overline{g_i^\leftarrow}, g_i^\leftarrow, f_i^\rightarrow) \simeq_d (f_i^\leftarrow,\overline{g_i^\leftarrow}, g_i^\leftarrow, g_i^\rightarrow) \xrightarrow{m(3,z_0+1)} (f_i^\leftarrow,\overline{g_i^\leftarrow}, f'_i), \]
	\[ (f_i^\leftarrow,\overline{g_i^\leftarrow}, f'_i) \simeq_{d} (f_i^\leftarrow,\overline{f_i^\leftarrow}, f'_i) \xrightarrow{d_{op(1)}} (f'_i). \]
	We conclude that $[f_i] = [f'_i]$, and thus $[F]=[f_1,\ldots,f'_i,\ldots,f_n]$.
\end{proof}

\begin{proposicion}
\label{prop:CutEqualTails}
	Let $X$ be a semi-coarse space and $F=(f_1,\ldots,f_n)$ be a string.
	Suppose that there exist $M,N\in \Z$ such that
	\[ f'_i(z) = f_i(z-M),\ f'_{i+1}(z) = f_{i+1}(z-N),\ f'_i\mid_{[0,\infty)} = \overline{f'_{i+1}}\mid_{[0,\infty)} , \]
	and there are no $M'<M$ or $N'>N$ which satisfies the same condition, for some $i\in\{1,\ldots,n-1\}$.
	Then we have that
	\[ [F] = [f_1, \ldots, f_{i-1}, g, f_{i+2},\ldots, f_n] \]
	with $g:\Z_1\to X$ defined by
	\begin{equation}
	\label{Eq:CutEqualTails}
	g(z)\coloneqq \begin{cases}
	f'_i(z) & z\geq 0,\\
	f'_{i+1}(z) & z<0
	\end{cases} 
	\end{equation}
\end{proposicion}

\begin{proof}
	Let $X$ and $F$ as in the statement.
	By \autoref{prop:Affinity}, we have that
	\[ [f_i,f_{i+1}] = [f'_i,f'_{i+1}]. \]
	Define the maps.
	\[ f_i^\leftarrow(z) \coloneqq \begin{cases}
	f'_i(z) & z\leq 0,\\
	f'_i(0) & z>0
	\end{cases}\hspace{5mm} f_i^\rightarrow(z) \coloneqq \begin{cases}
	f'_i(z) & z\geq 0,\\
	f'_i(0) & z<0
	\end{cases}\]
	then
	\[ f_{i+1}^\leftarrow(z) \coloneqq \begin{cases}
	f'_{i+1}(z) & z\leq 0,\\
	f'_{i+1}(0) & z>0
	\end{cases}\hspace{5mm} f_{i+1}^\rightarrow(z) \coloneqq \begin{cases}
	f'_{i+1}(z) & z\geq 0,\\
	f'_{i+1}(0) & z<0
	\end{cases}\]
	and, by hypothesis, we have that
	\[ f^\rightarrow_i = \overline{f^\leftarrow_{i+1}}. \]
	Thus
	\[ [f'_i,f'_{i+1}] = [f^\leftarrow_{i}, f^\rightarrow_{i}, f^\leftarrow_{i+1}, f^\rightarrow_{i+1}] = [f^\leftarrow_{i}, f^\rightarrow_{i+1}] = [g],\]
	where $g$ is defined as in \ref{Eq:CutEqualTails}. The result follows.
\end{proof}

	One residual effect of the finiteness of the values in the discrete intervals is the fundamental group's inability to distinguish some semi-coarse spaces which are clearly different;
	in particular, the semi-coarse fundamental group, and this equal for every semi-coarse homotopy group, is trivial for every coarse space.
	
	The following example illustrate two desirable properties: the fundamental groupoid is trivial for the indiscrete semi-coarse space and the fundamental groupoid is not necessarily trivial for coarse spaces.
	As a brief reminder, the indiscrete semi-coarse structure in $X$ is the collections $2^{X\times X}$.
	Although this is the structure of $X$ with more elements,
	it receives its name from a correspondence from topology, where the indiscrete topology is the structure which makes the identity set map continuous for every structure in the domain.

\begin{ejemplo}
\label{Exam:CoarFundGroupoidNotTriv}
	Let $\Z_{ind}$ be the integers joint with indiscrete semi-coarse structure.
	By \autoref{cor:ExtGroupoidSCInv}, we now that the fundamental groupoid of $\Z_{ind}$ only have one object.
	
	On the other hand, let $\Z_\infty$ be the canonical coarse structure in $\Z$, that is, $A$ is controlled in $\Z_\infty$ if there exists $R>0$ such that every pair $(x,y)\in A$ satisfies that  $d(x,y)\leq R$.
	Observe that $f,g:\Z_1\to \Z_\infty$ defined by $f(z) = |z|$ and $g(z)=2|z|$ are symmetric maps, but there is no a homotopy between them;
	we can note that if such a homotopy exists, we can add only a fixed finite number to the image of every point, but the distance between $|z|$ and $2|z|$ is always $|z|$, which is not fixed.
\end{ejemplo}

\begin{ejemplo}
	 For every $n\in\N$, $\pi_{\leq}(\Z_\infty)$ is not isomorphic to $\pi_{\leq}(\Z_n)$.
	 To prove it, let's define the following:	 
	 We say that an object $\langle f \rangle_\infty$ is \emph{eliminable} if whenever we have a string $F=(f_1,\ldots,f_n)$ ending in $\langle f \rangle_\infty$ and a string $G=(g_1,\ldots,g_m)$ beginning in $\langle f \rangle_\infty$, we can rewrite $F\star G$ as
	 \[H= (f_1,\ldots,f_{n-1},h,g_2,\ldots,g_m)\] 
	 through the \autoref{prop:CutEqualTails}.
	 We claim that the cardinality of eliminable objects is a semi-coarse invariant.
	 
	Observe that $\Z_n$ has one eliminable object, which is $\langle n|z| \rangle_\infty$;
	however, $\Z_\infty$ does not have any eliminable object.
	We can check that $\Z_\infty$ does not have eliminable objects taking any map $f:\Z_1\to \Z_\infty$ and defining $g$ as follows:
	\[ g\coloneq \begin{cases}
	f(-z)   & z\text{is an odd number}\\
	f(-z)-1 & z\text{is an even number}
	\end{cases} \]
	Then $(f)$ and $(g)$ satisfies the conditions, but we cannot apply \autoref{prop:CutEqualTails} to $(f,g)$.
	
	The existence of eliminable objects is also useful in distinguishing between two fundamental groupoids.
	If $\pi_\leq(X)$ has an eliminable object, $\langle f\rangle_\infty$, and $\pi_\leq(Y)$ has no eliminable objects, then some one-string $(g)$ ending in $\langle f\rangle_\infty$ and some one-string $(g')$ beginning in $\langle f\rangle_\infty$ are equivalent to some one-string $(h)$; however, this does not happen in $\pi_\leq(Y)$.
	Thus they are not isomorphic.
\end{ejemplo}

\begin{teorema}
\label{teo:FundGroupWithinFundGroupoid}
	Let $X$ be a semi-coarse space and $\pi_1^{sc}(X,x_0)$ the fundamental group at $x_0$ of $X$.
	Then there exists a monomorphism from $\pi_1^{sc}(X,x_0)$ to $\pi_{\leq 1}(X)$.
	In particular, if we see the group structure as a category with one object, we can send the unique object of $\pi_1^{sc}(X,x_0)$ to the object $\langle \hat{x}_0 \rangle_{\infty}$ with $\hat{x}_0:\Z\to X$ such that $\hat{x}_0(z)=x_0$.
\end{teorema}

\section{Relative Fundamental Groupoid and Van Kampen Theorem}
\label{cap.RelativeFundamentalGroupoid}

	The van Kampen Theorem is not true in general for the semi-coarse fundamental groupoid.
	To obtain a first impression of how problematic these maps can be, let's take $\Z_1$ well-split by the rays $(-\infty,0]$ and $[0,\infty)$.
	Next, we define the maps $g_n:\Z \to \Z$ such that $g_n(z)=n-|z-n|$ to build the following bornologous map
	\[g(z) \coloneqq \begin{cases}
		0 & z< 0\\
		(-1)^{n+1}g_n\left(z - 2\sum\limits_{i=0}^{n-1}i \right) & 2\sum\limits_{i=0}^{n-1}i \leq  z < 2\sum\limits_{i=0}^{n}i
	\end{cases}\]
	Look at \autoref{fig:NonVKT} to a representation of the map $g$. This map is a bornologous map which is not homotopy equivalent to a constant map.
	In addition, the points mapped through $g_i$ with an odd index are non-negative and with an even index are non-positive.
	This makes technically impossible to divide $g$ in a finite amount of strings whose images are either in $(-\infty,0]$ or $[0,\infty)$.
	
\begin{figure}[h!]
	\begin{tikzpicture}
		\foreach \x in {-3,-2,-1}
			{\filldraw[black] (\x/6,0) circle (0.2pt);}
		\foreach \x in {0,...,15}
			{\filldraw[Gray!30] (\x/2,0) circle (1.5pt) node[below,black] {\tiny{\inteval{\x-2}}};}
		\foreach \x in {0,1,2,4,8,14}
			{\filldraw[RedOrange] (\x/2,0) circle (1.5pt);}
		\filldraw[RedOrange] (3/2,1/2) circle (1.5pt);
		\filldraw[RedOrange] (9/2,1/2) circle (1.5pt);
		\filldraw[RedOrange] (13/2,1/2) circle (1.5pt);
		\filldraw[RedOrange] (5/2,-1/2) circle (1.5pt);
		\filldraw[RedOrange] (7/2,-1/2) circle (1.5pt);
		\filldraw[RedOrange] (15/2,-1/2) circle (1.5pt);
		\filldraw[RedOrange] (10/2,1) circle (1.5pt);
		\filldraw[RedOrange] (12/2,1) circle (1.5pt);
		\filldraw[RedOrange] (6/2,-1) circle (1.5pt);
		\filldraw[RedOrange] (11/2,3/2) circle (1.5pt);
		\foreach \x in {-3,-2,-1}
			{\filldraw[black] (\x/6+16/2+1/6,0) circle (0.2pt);}
\end{tikzpicture}
\caption{Representation of $g$}
\label{fig:NonVKT}
\end{figure}
	
	We can, however, restrict the tails of the maps from $\mathbb{Z}_1$ to get a relative version of the fundamental groupoid.
	We show in this section that a semi-coarse version of the van Kampen theorem does hold for the relative fundamental groupoid.
	To have good behavior, we need to ``control'' the tails of every function in a string. With this in mind, we present the following definition.	

\begin{definicion}
\label{def:TailsLivesIn}
	Let $X$ be a semi-coarse spaces and $\mathcal{U}\subset 2^X$.
	We say that the tails of $f:\Z \to X$ \emph{are controlled by} $\mathcal{U}$ if there exist $U_1,U_2\in \mathcal{U}$ such that the image of the left tail of $f$ is contained in $U_1$ and the image of the right tail of $f$ is contained in $U_2$.
	
	We say that a string $F=(f_1,\ldots,f_n)$ is controlled by $\mathcal{U}$ if the tails of every $f_i$ is controlled by $\mathcal{U}$.
\end{definicion}

	Observe we are not asking for $\mathcal{U}$ to have any specific properties; in particular,
	it does not have to be a cover or have any kind of intersection between its elements.
	
	Restricting $\pi_{\leq 1}(X)$ to the objects and morphisms of $\pi_{\leq 1}(X)$ for maps whose tails are controlled by $\mathcal{U}$ still results in a groupoid.
	This restriction does not affect the construction of the equivalence relation, only the set over what is defined;
	for example, $\mathcal{U}=\{ \{ x \} \mid x'\in X\}$ forces every map to have constant tails, but the relation in this class does not change at all.

\begin{definicion}
\label{def:URelativeFundamentalGroupoid}
	Let $X$ be a semi-coarse space and $\mathcal{U}\subset 2^X$.
	The \emph{$\mathcal{U}$-relative fundamental groupoid}, $\pi_{\leq 1}(X,\mathcal{U})$, is the category with the following objects and morphisms:
\begin{itemize}
	\item $Ob(\pi_{\leq 1}(X,\mathcal{U}))$ are the symmetric maps whose tails are controlled by $U\in\mathcal{U}$.
	
	\item $Hom_{\pi_{\leq 1}(X,\mathcal{U})}(f,g)$ are equivalence classes of strings $[F]=[f_1, \ldots, f_n]$ from $f$ to $g$ such that the tails of $F$ are controlled by $\mathcal{U}$.
\end{itemize}
\end{definicion}

	First, we introduce a lemma which we will need for the semi-coarse van Kampen theorem which mimics the implementation of the Lebesgue covering lemma \cite{Brown_2006}.

\begin{lema}
\label{lem:LebesgueCoveringSC}
	Let $X$ be a semi-coarse space, suppose that $A,B\subset X$ well-split $X$ and that $f:\mathbb{Z}_1\rightarrow X$ is bornologous.
	Additionally, suppose that each of the tails of $f$ are contained in one of $A$ or $B$.
	(The two tails may be contained in different sets.)
	Then, there exist bornologous functions $f_1,\ldots,f_n:\mathbb{Z}_1\rightarrow X$ such that  $f_i(\mathbb{Z})\subset A$ or $f_i(\mathbb{Z})\subset B$ for every $i\in 1,\dots,n$, $[f] = [f_1] \star [f_2] \star \ldots \star [f_n]$, and the right tail of $f_i$ is the same as the left tail of $f_{i+1}$.
\end{lema}

\begin{proof}
	Let $X,A,B,$ and $f$ be as in the statement of the lemma.
	We treat the case where the left tail of $f$ is eventually in $A$ and the right tail of $f$ is eventually in $B$. The other cases are handled similarly.

	If $f(\mathbb{Z})\subset A$ or $f(\mathbb{Z})\subset B$, then the conclusions of the lemma are satisfied.
	Suppose now that $f(\mathbb{Z}) \not \subset A$ and $f(\mathbb{Z}) \not \subset B$.
	Then there exists a minimum $z\in\mathbb{Z}$ such that $f(z)\in B-A$.
	We repeat this searching to find the set $z_1<z_2<\ldots<z_{k}$ such that $f(z_i)\in B-A$ and $f(z_{i+1})\in A-B$ for every $i\in \{1,3,\ldots,2\lfloor \frac{k}{2} \rfloor +1\}$.

\begin{figure}[h]
\begin{tikzpicture}
\foreach \x in {-3,...,-1}
	{\filldraw[black] (\x/6,0) circle (0.2pt);}
\foreach \x in {0,...,5}
	{\filldraw[black] (\x/3,0) circle (1.5pt);}
\foreach \x in {0,...,4}
	{\filldraw[blue] (\x/3+6/3,0) circle (1.5pt);}
\filldraw[red] (11/3,0) circle (1.5pt);
\foreach \x in {0,...,5}
	{\filldraw[Green] (\x/3+12/3,0) circle (1.5pt);}
\filldraw[Purple] (18/3,0) circle (1.5pt);
\foreach \x in {0,...,3}
	{\filldraw[blue] (\x/3+19/3,0) circle (1.5pt);}
\foreach \x in {0,...,2}
	{\filldraw[black] (\x/6+22/3+1/6,0) circle (0.2pt);}
\foreach \x in {0,...,2}
	{\filldraw[Green] (\x/3+22/3+4/6,0) circle (1.5pt);}
\foreach \x in {0,...,3}
	{\filldraw[black] (\x/3+25/3+4/6,0) circle (1.5pt);}
\foreach \x in {0,...,2}
	{\filldraw[black] (\x/6+28/3+5/6,0) circle (0.2pt);}
	
\filldraw[black] (3/3-1/6,.5+2/6) circle (0pt) node[anchor=south]{\footnotesize{$A$}};
\filldraw[black] (3/3-1/6,.25+2/6) circle (0pt) node[anchor=south]{\footnotesize{Left tail}};
\foreach \x in {-3,...,-1}
	{\filldraw[black] (\x/6,.25+1/6) circle (0.2pt);}
\draw[thick,black] (0,.25+1/6) -- (3/3-2/6,.25+1/6);
\draw[thick,black] (3/3-2/6,.25+1/6) arc (270:360:1/6);
\draw[thick,black] (3/3-1/6,.25+2/6) arc (180:270:1/6);
\draw[thick,black] (3/3,.25+1/6) -- (5/3,.25+1/6);
\draw[thick,black] (5/3+1/6,.25) arc (0:90:1/6);

\filldraw[black] (8/3,.25+2/6) circle (0pt) node[anchor=south]{\footnotesize{$A$}};
\draw[thick,blue] (6/3,.25+1/6) arc (90:180:1/6);
\draw[thick,blue] (6/3,.25+1/6) -- (7/3+1/6,.25+1/6);
\draw[thick,blue] (7/3+1/6,.25+1/6) arc (270:360:1/6);
\draw[thick,blue] (8/3,.25+2/6) arc (180:270:1/6);
\draw[thick,blue] (8/3+1/6,.25+1/6) -- (10/3,.25+1/6);
\draw[thick,blue] (10/3+1/6,.25) arc (0:90:1/6);

\filldraw[black] (11/3,-.25-2/6) circle (0pt) node[anchor=north]{\footnotesize{$f(z_1)$}};
\draw[thick,red] (11/3,-.1) -- (11/3,-.25-2/6);
\draw[thick,red] (11/3,-.1) -- (11/3+.05,-.1-.1);
\draw[thick,red] (11/3,-.1) -- (11/3-.05,-.1-.1);

\filldraw[black] (14/3,.25+2/6) circle (0pt) node[anchor=south]{\footnotesize{$B$}};
\draw[thick,Green] (11/3,.25+1/6) arc (90:180:1/6);
\draw[thick,Green] (11/3,.25+1/6) -- (14/3-1/6,.25+1/6);
\draw[thick,Green] (14/3-1/6,.25+1/6) arc (270:360:1/6);
\draw[thick,Green] (14/3,.25+2/6) arc (180:270:1/6);
\draw[thick,Green] (15/3-1/6,.25+1/6) -- (17/3,.25+1/6);
\draw[thick,Green] (17/3+1/6,.25) arc (0:90:1/6);

\filldraw[black] (18/3,-.25-2/6) circle (0pt) node[anchor=north]{\footnotesize{$f(z_2)$}};
\draw[thick,Purple] (18/3,-.1) -- (18/3,-.25-2/6);
\draw[thick,Purple] (18/3,-.1) -- (18/3+.05,-.1-.1);
\draw[thick,Purple] (18/3,-.1) -- (18/3-.05,-.1-.1);

\filldraw[black] (28/3+1/6,.5+2/6) circle (0pt) node[anchor=south]{\footnotesize{$B$}};
\filldraw[black] (28/3+1/6,.25+2/6) circle (0pt) node[anchor=south]{\footnotesize{Right tail}};
\draw[thick,black] (27/3,.25+1/6) arc (90:180:1/6);
\draw[thick,black] (27/3,.25+1/6) -- (28/3,.25+1/6);
\draw[thick,black] (28/3,.25+1/6) arc (270:360:1/6);
\draw[thick,black] (28/3+1/6,.25+2/6) arc (180:270:1/6);
\draw[thick,black] (29/3,.25+1/6) -- (30/3,.25+1/6);
\foreach \x in {0,...,2}
	{\filldraw[black] (\x/6+28/3+5/6,.25+1/6) circle (0.2pt);}
\end{tikzpicture}
\caption{Dividing $f$. Black points represent the tails of $f$, blue points are in $A$, and green points are in $B$. The red points represent that $f(z_i)\in B-A$ and the purple points represent that $f(z_i)\in A-B$.}
\label{fig:DividingFprime}
\end{figure}

	We are going to define recursively an $\simeq$-equivalent map to $f$, say $f'$, which satisfies the conclusions of the lemma.
		
	For each $i\in \{1,\ldots, k\}$ we have two cases: Whether $f(z_i-1)\in A\cap B$ or not.
	\begin{itemize}
		\item If $f(z_i-1)\in A\cap B$, define $\lambda_i:\{0,1,2\}\rightarrow X$ as $\lambda_i(0)=f(z_i-1)$, $\lambda_i(1)=f(z_i-1)$ and $\lambda_i(2)=f(z_i)$. Observe that its image is in $A\sqcup_{A\cap B}B$.
		
		\item If $f(z_i-1)\notin A\cap B$. Since $A,B$ well-split $X$, there exists a short path $\lambda_i:\{0,1,2\}\rightarrow X$ whose image is in $A\sqcup_{A\cap B}B$ and is such that $\lambda_i(0)=f(z_i-1)$, $\lambda_i(2)=f(z_i)$ and $\lambda_i(1)\in A\cap B$ (the last contention is for \autoref{prop:PathBetweenElements}).
	\end{itemize}
		
	We define $f'_i:\Z_1 \to X$ as follows
	\begin{align*}
	f'_i(z) \coloneqq \begin{cases}
		f'_{i-1}(z+1)  &  z \leq z_i-2\\
		\lambda_i(1) &  z = z_i-1\\
		f'_{i-1}(z)    &  z \geq z_i
	\end{cases}
	\end{align*}
	for every $i\in\{1,\ldots,k\}$ and with $f'_0\coloneqq f$ and $f'\coloneqq f'_k$. We note that in every step, we are just adding $\lambda_i(1)$ between $f(z_i-1)$ and $f(z_i)$.
	
	For each $i = 1,\dots,n$, define $z'_i\coloneqq z_i-1-(k-i)$.	
	Hence we obtain a new sequence $z'_1<z'_2<\ldots<z'_k$ such that $f(z_i')\in A\cap B$, and $f(z'_i+1)\in B-A$ and $f(z'_{i+1}+1)\in A-B$ for every $i\in\{1, 3,\ldots, 2\lfloor\frac{k}{2}\rfloor+1 \}$.

\begin{figure}[h]
\begin{tikzpicture}
\foreach \x in {-3,...,-1}
	{\filldraw[black] (\x/6,0) circle (0.2pt);}
\foreach \x in {0,...,5}
	{\filldraw[black] (\x/3,0) circle (1.5pt);}
\foreach \x in {0,...,4}
	{\filldraw[blue] (\x/3+6/3,0) circle (1.5pt);}
\filldraw[Purple] (11/3,0) circle (1.5pt);
\foreach \x in {0,...,5}
	{\filldraw[Green] (\x/3+12/3,0) circle (1.5pt);}
\filldraw[Purple] (18/3,0) circle (1.5pt);
\foreach \x in {0,...,3}
	{\filldraw[blue] (\x/3+19/3,0) circle (1.5pt);}
\foreach \x in {0,...,2}
	{\filldraw[black] (\x/6+22/3+1/6,0) circle (0.2pt);}
\foreach \x in {0,...,2}
	{\filldraw[Green] (\x/3+22/3+4/6,0) circle (1.5pt);}
\foreach \x in {0,...,3}
	{\filldraw[black] (\x/3+25/3+4/6,0) circle (1.5pt);}
\foreach \x in {0,...,2}
	{\filldraw[black] (\x/6+28/3+5/6,0) circle (0.2pt);}
	
\filldraw[black] (3/3-1/6,.5+2/6) circle (0pt) node[anchor=south]{\footnotesize{$U$}};
\filldraw[black] (3/3-1/6,.25+2/6) circle (0pt) node[anchor=south]{\footnotesize{Left tail}};
\foreach \x in {-3,...,-1}
	{\filldraw[black] (\x/6,.25+1/6) circle (0.2pt);}
\draw[thick,black] (0,.25+1/6) -- (3/3-2/6,.25+1/6);
\draw[thick,black] (3/3-2/6,.25+1/6) arc (270:360:1/6);
\draw[thick,black] (3/3-1/6,.25+2/6) arc (180:270:1/6);
\draw[thick,black] (3/3,.25+1/6) -- (5/3,.25+1/6);
\draw[thick,black] (5/3+1/6,.25) arc (0:90:1/6);

\filldraw[black] (8/3,.25+2/6) circle (0pt) node[anchor=south]{\footnotesize{$A$}};
\draw[thick,blue] (6/3-1/3,.25+1/6) -- (7/3+1/6,.25+1/6);
\draw[thick,blue] (7/3+1/6,.25+1/6) arc (270:360:1/6);
\draw[thick,blue] (8/3,.25+2/6) arc (180:270:1/6);
\draw[thick,blue] (8/3+1/6,.25+1/6) -- (11/3,.25+1/6);
\draw[thick,blue] (11/3+1/6,.25) arc (0:90:1/6);

\filldraw[black] (11/3,-.25-2/6) circle (0pt) node[anchor=north]{\footnotesize{$f'(z_1)$}};
\draw[thick,Purple] (11/3,-.1) -- (11/3,-.25-2/6);
\draw[thick,Purple] (11/3,-.1) -- (11/3+.05,-.1-.1);
\draw[thick,Purple] (11/3,-.1) -- (11/3-.05,-.1-.1);

\filldraw[black] (14/3+1/6,.25+2/6) circle (0pt) node[anchor=south]{\footnotesize{$B$}};
\draw[thick,Green] (11/3,.25+1/6) arc (90:180:1/6);
\draw[thick,Green] (11/3,.25+1/6) -- (14/3,.25+1/6);
\draw[thick,Green] (14/3,.25+1/6) arc (270:360:1/6);
\draw[thick,Green] (14/3+1/6,.25+2/6) arc (180:270:1/6);
\draw[thick,Green] (15/3,.25+1/6) -- (18/3,.25+1/6);
\draw[thick,Green] (18/3+1/6,.25) arc (0:90:1/6);

\filldraw[black] (18/3,-.25-2/6) circle (0pt) node[anchor=north]{\footnotesize{$f'(z_2)$}};
\draw[thick,Purple] (18/3,-.1) -- (18/3,-.25-2/6);
\draw[thick,Purple] (18/3,-.1) -- (18/3+.05,-.1-.1);
\draw[thick,Purple] (18/3,-.1) -- (18/3-.05,-.1-.1);

\filldraw[black] (28/3+1/6,.5+2/6) circle (0pt) node[anchor=south]{\footnotesize{$V$}};
\filldraw[black] (28/3+1/6,.25+2/6) circle (0pt) node[anchor=south]{\footnotesize{Right tail}};
\draw[thick,black] (27/3,.25+1/6) arc (90:180:1/6);
\draw[thick,black] (27/3,.25+1/6) -- (28/3,.25+1/6);
\draw[thick,black] (28/3,.25+1/6) arc (270:360:1/6);
\draw[thick,black] (28/3+1/6,.25+2/6) arc (180:270:1/6);
\draw[thick,black] (29/3,.25+1/6) -- (30/3,.25+1/6);
\foreach \x in {0,...,2}
	{\filldraw[black] (\x/6+28/3+5/6,.25+1/6) circle (0.2pt);}
\end{tikzpicture}
\caption{Dividing $f'$}
\label{fig:DividingFi}
\end{figure}

	We now define the following maps for $i\in\{1,\ldots k-1\}$
	\begin{align*}
		f_{i}(z)\coloneqq\left\lbrace\begin{array}{ll}
			f'(z'_i)     & \text{ if } z<z'_i\\
			f'(z)        & \text{ if } z'_i\leq z \leq z'_{i+1}\\
			f'(z'_{i+1}) & \text{ if } z>z'_{i+1} .
		\end{array}\right.
	\end{align*}
Finally, we define the maps $f_L$ and $f_R$ by
\begin{align*}
f_{L}(z)\coloneqq\left\lbrace\begin{array}{ll}
f'(z)    & \text{ if } z\leq z'_1\\
f'(z'_1) & \text{ if } z>z'_1,
\end{array}\right. &
f_{R}(z)\coloneqq\left\lbrace\begin{array}{ll}
f'(z'_k) & \text{ if } z<z'_k\\
f'(z)    & \text{ if } z\geq z'_k,
\end{array}\right.
\end{align*}

	Hence, by construction $f_L(\mathbb{Z})\subset A$, $f_i(\mathbb{Z})\subset B$ and $f_{i+1}(\mathbb{Z})\subset A$ for $i\in \{1, 3,\ldots\}$, and $f_{R} \subset B$.
	In addition, the right tail of $f_i$ and the left tail of $f_{i+1}$ are the same and constant.
	
	Defining the $(k+1)$-string $F$ by $F\coloneqq (f_L, f_1, f_2, \ldots, f_{k-1}, f_R)$ and we observe that $F\simeq (f)$ by merging the constant tails of every map in $F$.
\end{proof}

	Before the proof of the main theorem, we prove the van Kampen theorem for a connected semi-coarse space; this implies in particular that $A\cap B \neq \varnothing$ if $A,B$ well-split $X$, since \autoref{prop:SplitWell_IFF_PathConnected}.

\begin{teorema}
\label{theo:VanKampenConnected}
	Let $X$ be a connected semi-coarse space, and let $\mathcal{U}\subset 2^X$ be such that $U\in\mathcal{U}$ is connected.
	Consider $A,B\subset X$ which satisfies
\begin{enumerate}
	\item $A,B$ well-split $X$,
	
	\item for every connected component $V$ of $A$, $B$ or $A\cap B$, there exists a set $U \in \mathcal{U}$ such that $U\cap V \neq \varnothing$, and
	
	\item for every $U\in\mathcal{U}$, $U\subset A$ or $U\subset B$.
\end{enumerate}
	Then the diagram
\begin{align*}
\xymatrix{
\pi_{\leq 1}(A\cap B, \mathcal{U}\cap A\cap B) \ar[r]^{i_A} \ar[d]^{i_B} &  \pi_{\leq 1}(A, \mathcal{U}\cap A) \ar[d]^{i_1}\\
\pi_{\leq 1}(B, \mathcal{U}\cap B) \ar[r]^{i_2} & \pi_{\leq 1}(X, \mathcal{U})
}
\end{align*}
is a pushout.
\end{teorema}

\begin{proof}
	Let $X$ be a connected semi-coarse space and $\mathcal{U}\subset 2^X$such that $U\in\mathcal{U}$ is connected. Let $G$ be a groupoid such that the following diagram commutates:

\begin{figure}[h]
\[
\xymatrix{
\pi_{\leq 1}(A\cap B, \mathcal{U}\cap A\cap B) \ar[r]^{i_A} \ar[d]^{i_B} &  \pi_{\leq 1}(A, \mathcal{U}\cap A) \ar[d]^{i_1} \ar@/^1pc/[ddr]^{j_1}\\
\pi_{\leq 1}(B, \mathcal{U}\cap B) \ar[r]^{i_2} \ar@/_1pc/[rrd]^{j_2} & \pi_{\leq 1}(X,\mathcal{U}) \ar@{..>}[rd]^{h} \\
 & & G
}
\]
\caption{}\label{diag:VKConnPO}
\end{figure}

	We proceed to construct a map (functor) $h:\pi_{\leq 1}(X,\mathcal{U})\rightarrow G$ such that is the unique such that $j_2=hi_2$ and $j_1=hi_1$.

	We first work with the objects. Let $f:\Z_1 \to X$ be a symmetric map whose tails are controlled by $\mathcal{U}$. Then there exists $U\in\mathcal{U}$ such that the tails of $f$ are contained in $U$, and therefore also in $A$ or $B$.
	Without loss of generality, we can say that $f(\Z)\subset A$ or $f(\Z)\subset B$; this is because there exists $f'\in \langle f\rangle_\infty$ such the image is within one of these two sets.
	Therefore, for every $f:\Z_1 \to X$, we may assume that either $\langle f \rangle_\infty \in Ob(\pi_{\leq 1}(A,\mathcal{U}\cap A)$ or $\langle f \rangle_\infty \in Ob(\pi_{\leq 1}(B,\mathcal{U}\cap B)$.
	
	Define $h:Ob(\pi_{\leq 1}(X,\mathcal{U}))\rightarrow Ob(G)$ such that $h(\langle f \rangle_\infty)=j_1(\langle f \rangle_\infty)$ if $\langle f \rangle_\infty \in Ob(\pi_{\leq 1}(A,\mathcal{U}\cap A))$ and $h(\langle f \rangle_\infty)=j_2(\langle f \rangle_\infty)$ if $\langle f \rangle_\infty\in Ob(\pi_{\leq 1}(B,\mathcal{U}\cap B))$.
	In case the tails of $f$ are contained in $A\cap B$, the map is well-defined because $j_2 i_B = j_1 i_A$.
	This map is unique in the objects by construction.

	We now consider the morphisms. Let $F=(f_1,\ldots,f_n)$ be a string. By \autoref{lem:LebesgueCoveringSC}, we can write
	\[[F] = [f_{1,L}]\star [f'_{1,1}]\star \ldots\star [f'_{1,k_1}]\star [f_{1,R}]\star [f_{2,L}]\star [f'_{2,1}]\star\ldots\star[f_{n,R}]\]
	such that each tail of each factor is contained in $A$ or is contained in $B$.
	We pursue to build factors which each one is a morphism in $\pi_{\leq 1}(A,\mathcal{U}\cap A)$ or in $\pi_{\leq 1}(B,\mathcal{U}\cap B)$.
	By the construction in \autoref{lem:LebesgueCoveringSC}, every element of this product has constant tails, except the left tail of $f_{i,L}$ and the right tail of $f_{i,R}$ for every $i$.
	In addition, since, by construction, these constant tails are in a
	connected component $V$ of $A\cap B$, there exists by hypothesis a set $U_{i,j}\in\mathcal{U}$ which meets $V$ non-trivially.

\begin{figure}[h!]
\begin{tikzpicture}
    \filldraw [color=Red, opacity=0.6] plot[smooth cycle] coordinates {(0, 0) (2, 0.5) (4, .25) (5, 2) (4,3.5) (3,3) (2,3.5) (1,2.5) (-.5,2)};
    \filldraw [color=blue, opacity = 0.60] plot[smooth cycle] coordinates {(2, 1) (4, 1) (4.5, 2) (4, 3) (3,2.5) (2,2.5) (1,1.5)};
    \filldraw [color=blue, opacity = 0.60] plot[smooth cycle] coordinates {(6,2) (6.5,1) (7,2.25) (6.3,2.8)};
    
    \filldraw[Black] (0.3,1.4) circle (1.5pt);
    \filldraw[Gray!50] (0.5,1.42) circle (1pt);
    \filldraw[Gray!50] (0.7,1.43) circle (1pt);
    \filldraw[Gray!50] (0.9,1.47) circle (1pt);
    \filldraw[Gray!50] (1.1,1.55) circle (1pt);
    \filldraw[Gray!50] (1.3,1.64) circle (1pt);
    \filldraw[Gray!50] (1.5,1.73) circle (1pt);
    \filldraw[Gray!50] (1.7,1.85) circle (1pt);
    \filldraw[Cyan] (1.9,2) circle (1.5pt);
    
	\filldraw[White] (0.1,-0.3) rectangle (2.1,-0.7);   
	
	\draw[->, very thick] ( 1.1,-0.3) -- (1.1,1.4);
    
    \filldraw[Black] (0.3,-0.5) circle (1.5pt);
    \filldraw[Gray!50] (0.5,-0.5) circle (1pt);
    \filldraw[Gray!50] (0.7,-0.5) circle (1pt);
    \filldraw[Gray!50] (0.9,-0.5) circle (1pt);
    \filldraw[Gray!50] (1.1,-0.5) circle (1pt);
    \filldraw[Gray!50] (1.3,-0.5) circle (1pt);
    \filldraw[Gray!50] (1.5,-0.5) circle (1pt);
    \filldraw[Gray!50] (1.7,-0.5) circle (1pt);
    \filldraw[Cyan] (1.9,-0.5) circle (1.5pt);

    \filldraw[black] (0.3,-0.7) circle (0pt) node[anchor=north]{$0$};
    \filldraw[black] (1.9,-0.7) circle (0pt) node[anchor=north]{$m_{i,j}$};
    \filldraw[black] (1.1,0.55) circle (0pt) node[anchor=east]{$\gamma'_{i,j}$};
    
    \filldraw[Red, opacity = 0.6] (5,-0.5) circle (3pt) node[anchor=west, black, opacity = 1]{$V_{}$};
    \filldraw[blue, opacity = 0.6] (6,-0.5) circle (3pt) node[anchor=west, black, opacity = 1]{$U_{i,j}$};
    \filldraw[Red, opacity = 0.6] (7,-0.5) circle (3pt) node[anchor=west, black, opacity = 1]{$V\cap U_{i,j}$};\filldraw[blue, opacity = 0.6] (7,-0.5) circle (3pt);
\end{tikzpicture}
\caption{The black point represents one constant tail of some $f_{i,j}$. The red stain represents $V$, the connected component of $A\cap B$, the blue stain represents $U_{i,j}$ and the purple stain represent the intersection of both. The yellow path is one possible path from the tail to $U_{i,j}$}
\label{fig:GammaPathComponent}
\end{figure}

	Then, since $V \subset A\cap B$ is connected, there exists a bornologous map $\gamma'_{i,j}:\{0,1,\ldots,m_{i,j}\}\rightarrow X$ whose image is totally contained in $V$, such that  $\gamma'_{i,j}(0)$ takes the constant value of the right tail of $f_{i,j}$ and $f_{i,j}(m_{i,j})\in U_{i,j}$.
	We extend $\gamma'_{i,j}$ to define it over all of the integers by 
\begin{align*}
\gamma_{i,j}(z)\coloneqq \left\lbrace \begin{array}{ll}
\gamma'_{i,j}(0) & \text{ if }z\leq 0\\
\gamma'_{i,j}(z) & \text{ if }0 < z < m_{i,j}\\
\gamma'_{i,j}(m_{i,j}) & \text{ if }z\leq m_{i,j}.
\end{array}\right.
\end{align*}
	By convention, we denote $f_{i,0}\coloneqq f_{i,L}$ and $f_{i,k_n+1}\coloneqq f_{i,R}$.
	Therefore we take
	\[ [g_{i,j}] \coloneqq [\gamma_{i,j-1}] \star [f_{i,j}] \star [\overline{\gamma_{i,j}}] \text{ for } (i,j)\notin \{ (i,0), (i,k_i+1) \mid i\in \{1,\ldots,n\} \},\]
	and we define $[g_{i,0}]\coloneqq [g_{i,L}]\coloneqq [f_{i,L}]\star [\gamma_{i,0}]$ and $[g_{i,k_i+1}]\coloneqq [g_{i,R}]\coloneqq[\overline{\gamma_{i,k_i}}]\star[f_{i,R}]$ for every $i\in\{1,\ldots, n\}$.

	Observe that the left tail of each $f_{i,L}$ and the right tail of each $f_{i,R}$ is already contained in $A$ or $B$.
	
	Every $g_{i,j}$ is constructed to be completely contained in $A$ or $B$, so $[g_{i,j}]$ is also a string in $A$ or in $B$.
	In addition, each of its tails is contained in some $U\in\mathcal{U}$; therefore $g_{i,j}\in \pi_{\leq 1}(A,\mathcal{U}\cap A)$ or $g_{i,j}\in \pi_{\leq 1}(B,\mathcal{U}\cap B)$.
	
	For the diagram in \autoref{diag:VKConnPO} to commute, the morphism $h$ must satisfy 
	\begin{align}\label{Eq:DefinitionH}
h[F] = \kappa_{1,0}[g_{1,L}]\star \kappa_{1,1}[g_{1,1}]\star\ldots \kappa_{n,k_n+1}[g_{n,R}]
\end{align}
where $\kappa_{i,j} = j_1$ if the image of $g_{i,j}$ is in $A$ and 
$\kappa_{i,j} = j_2$ if the image of $g_{i,j}$ is in $B$.
	Defining $h$ in this way would assure 
the uniqueness of $h$ and make $\pi_{\leq 1}(X,\mathcal{U})$ the pushout.
	The remainder of this proof will be to show that this construction is well-defined, i.e. that $h$ does not depend on the
	choice of representative of $[F]$ or on the choices of the functions $\gamma'_{i,j}$.
	
	To prove that the construction does not depend of the election of the $\gamma'_{i,j}$, we consider a $U'_{i,j}$ which meets the component $V$ and a function $\eta'_{i,j}:\{0,\ldots,l_{i,j}\}\rightarrow X$ whose image is totally contained in $V$ and such that 
	\begin{enumerate}[label=(\alph*)]
		\item $\eta'_{i,j}(0)$ is the value taken by the right tail of $f_{i,j}$, and
		\item $\eta'_{i,j}(l_{i,j})\in U'_{i,j}$. 
	\end{enumerate}
	We define $\eta_{i,j}$ in the same way we defined $\gamma_{i,j}$ using $\eta'_{i,j}$ in place of $\gamma'_{i,j}$, and we observe that:
\begin{align*}
	\kappa_{i,j}[g_{i,j}]\star\kappa_{i^*,j^*}[g_{i^*,j^*}] = & \kappa_{i,j}[\overline{\gamma_{i_*,j_*}},f_{i,j},\gamma_{i,j},\overline{\gamma_{i,j}},\eta_{i,j},\overline{\eta_{i,j}},\gamma_{i,j}]\star \kappa_{i^*,j^*}[\overline{\gamma_{i,j}},f_{i^*,j^*}\gamma_{i^*,j^*}]\\
	= & \kappa_{i,j}[\overline{\gamma_{i_*,j_*}},f_{i,j},\eta_{i,j},\overline{\eta_{i,j}},\gamma_{i,j}] \star \kappa_{i^*,j^*}[\overline{\gamma_{i,j}},f_{i^*,j^*}\gamma_{i^*,j^*}]\\
	= & \kappa_{i,j}[\overline{\gamma_{i_*,j_*}},f_{i,j},\eta_{i,j}] \star \kappa_{i^*,j^*}[\overline{\eta_{i,j}},\gamma_{i,j},\overline{\gamma_{i,j}},f_{i^*,j^*}\gamma_{i^*,j^*}]\\
	= & \kappa_{i,j}[\overline{\gamma_{i_*,j_*}},f_{i,j},\eta_{i,j}] \star \kappa_{i^*,j^*}[\overline{\eta_{i,j}},f_{i^*,j^*}\gamma_{i^*,j^*}].
\end{align*}

	with $(i^*,j^*)$ and $(i_*,j_*)$ the immediately subsequent and preceding pairs, respectively, of $(i,j)$ in the lexicographic order.

	With a similar argument, we observe that the image neither depends on the election of the value of $f'(z'_k)$ in the construction \autoref{lem:LebesgueCoveringSC};
	we use part of the notation from that lemma in this fragment of the proof.
	Consider another function $\lambda'_{k}:\{0,1,2\}\rightarrow X$ in $A\sqcup_{A\cap B} B$ such that $\lambda'_k(0)=f(z_k-1)$ and $\lambda'_k(2)=f(z_k)$.
	Since $A,B$ well-split $X$, there exists a path $\alpha:\{0,\ldots,m'\}\rightarrow A\cap B$ such that $\alpha(0)=f'(z'_k)$ and $\alpha(m')=\lambda'_k(1)$.
	We define $\eta'_{i,j}\coloneqq \overline{\alpha} \star \gamma'_{i,j}$, $f''_{i,j}$ equal to $f_{i,j}$ but replacing the value of its right tail by $\lambda'_k(1)$, and we set $f''_{i^*,j^*}$ equal to $f_{i^*,j^*}$ but replacing the value of its left tail by $\lambda'_k(1)$.
	Thus
\begin{align*}
	\kappa_{i,j}[g_{i,j}] \star \kappa_{i^*,j^*}[g_{i^*,j^*}] = &  \kappa_{i,j}[\overline{\gamma_{i_*,j_*}}, f_{i,j},\gamma_{i,j}] \star \kappa_{i^*,j^*}[\overline{\gamma_{i,j}}, f_{i^*,j^*},\gamma_{i^*,j^*}]\\
	 = & \kappa_{i,j}[\overline{\gamma_{i_*,j_*}}, f''_{i,j},\eta_{i,j}] \star \kappa_{i^*,j^*}[\overline{\eta_{i,j}}, f''_{i^*,j^*},\gamma_{i^*,j^*}]
\end{align*}
The last equality follows from the graph \autoref{fig:GraphDeletingAlpha}, which is given by the properties of splitting-well.

	The rest of the $\gamma_{i,j}$ are also replaced  with $\eta_{i,j}$ in the whole expression of the factor of $[F]$.
	We only write that part of the expression as an example.

\begin{figure}[h]
\begin{tikzpicture}
\filldraw[black] (0,0) circle (0pt) node[anchor=south]{$\alpha(i)$};
\filldraw[black] (0,-2) circle (0pt) node[anchor=north]{$\alpha(i+1)$};
\filldraw[black] (-1.732,-1) circle (0pt) node[anchor=east]{$\lambda'_k(0)$};
\filldraw[black] (1.732,-1) circle (0pt) node[anchor=west]{$\lambda'_k(2)$};

\draw[RoyalBlue,line width=0.5mm] (-1.732,-1) -- (0,0);
\draw[RoyalBlue,line width=0.5mm] (-1.732,-1) -- (0,-2);
\draw[Red,line width=0.5mm] (1.732,-1) -- (0,0);
\draw[Red,line width=0.5mm] (1.732,-1) -- (0,-2);
\draw[Orchid,line width=0.5mm] (0,0) -- (0,-2);

\filldraw[Orchid] (0,0) circle (2pt);
\filldraw[Orchid] (0,-2) circle (2pt);
\filldraw[RoyalBlue] (-1.732,-1) circle (2pt);
\filldraw[Red] (1.732,-1) circle (2pt);
\end{tikzpicture}
\caption{Deleting $\alpha$'s}
\label{fig:GraphDeletingAlpha}
\end{figure}

	We claim that $h$ does not depend on the representative of the class $[F]$.
	To see it, it is enough to work with each of the equivalence classes which define $\simeq$ in \autoref{def:FundamentalGroupoid}.
	We check case by case:
	
	The equivalence relation given by merging tails ($\xleftrightarrow{m(i)}$) which are constant and equal to each other also does not present a problem in this construction:
	We only have to observe that tails are in the same set and how they combine.
	That is, let $[f,g]$ such that we can merge their related tails and call $[f']$ the result of merging these two functions.
	The case when $f$ or $g$ are constant functions is trivial, so let's assume that $f$ and $g$ are not constant functions.
	Then there exist $\alpha,\beta\in \Z$ such that $f(\alpha)\neq f(\alpha-1)$, $f(i)=f(i+1)$ for every $k\geq \alpha$, $g(\beta)\neq g(\beta+1)$ and $g(k-1)=g(k)$ for every $k\leq\beta$.
	Consider $z_1,z_2,\ldots,z_{n_f}$ and $y_1,y_2,\ldots, y_{n_g}$ the points given by \autoref{lem:LebesgueCoveringSC} for $f$ and $g$, respectively.
	By construction, $z_{n_f}\leq \alpha$ and $y_1> \beta$.
	This implies that $f'$ has the points $w_1,\ldots w_{n_f+n_g}$ given by \autoref{lem:LebesgueCoveringSC} represented in the \autoref{fig:MergingVanKampen}.
	\begin{figure}[h]
	\begin{tikzpicture}
		\filldraw[RoyalBlue] (0,0) circle (2pt);
		\draw[->] (0,0.5) node[anchor=south,black]{$z_1$} -- (0,0.1);
		\filldraw[RoyalBlue] (1/2,0) circle (2pt);
		\draw[->] (1/2,0.5) node[anchor=south,black]{$z_2$} -- (1/2,0.1);
		\filldraw[Black] (1/2+1/6,0) circle (0.5pt);
		\filldraw[Black] (1/2+2/6,0) circle (0.5pt);
		\filldraw[Black] (1/2+3/6,0) circle (0.5pt);
		\filldraw[RoyalBlue] (1+1/6,0) circle (2pt);
		\draw[->] (1+1/6,0.5) node[anchor=south,black]{$z_{n_f}$} -- (1+1/6,0.1);
		
		\filldraw[Red] (2,0) circle (2pt);
		\draw[->] (2,0.5) node[anchor=south,black]{$y_1$} -- (2,.1);
		\filldraw[Red] (5/2,0) circle (2pt);
		\draw[->] (5/2,0.5) node[anchor=south,black]{$y_2$} -- (5/2,.1);
		\filldraw[Black] (5/2+1/6,0) circle (0.5pt);
		\filldraw[Black] (5/2+2/6,0) circle (0.5pt);
		\filldraw[Black] (5/2+3/6,0) circle (0.5pt);
		\filldraw[Red] (3+1/6,0) circle (2pt);
		\draw[->] (3+1/6,0.5) node[anchor=south,black]{$y_{n_g}$} -- (3+1/6,0.1);
		
		\draw[->,very thick] (0,-.1) -- (0+1/6,-0.9);
		\draw[->,very thick] (1/2,-.1) -- (1/2+1/6,-0.9);
		\draw[->,very thick] (1+1/6,-.1) -- (1+1/6+1/6,-0.9);
		\draw[->,very thick] (2,-.1) -- (1+1/6+1/2+1/6,-0.9);
		\draw[->,very thick] (5/2,-.1) -- (1+1/6+1+1/6,-0.9);
		\draw[->,very thick] (3+1/6,-.1) -- (1+1/6+3/2+1/6+1/6,-0.9);
		
		\filldraw[Orchid] (0+1/6,-1) circle (2pt);
		\draw[->] (0+1/6,-1.5) node[anchor=north,black]{$w_1$} -- (0+1/6,-1.1);
		\filldraw[Orchid] (1/2+1/6,-1) circle (2pt);
		\draw[->] (1/2+1/6,-1.5) node[anchor=north,black]{$w_2$} -- (1/2+1/6,-1.1);
		\filldraw[Black] (1/2+1/6+1/6,-1) circle (0.5pt);
		\filldraw[Black] (1/2+2/6+1/6,-1) circle (0.5pt);
		\filldraw[Black] (1/2+3/6+1/6,-1) circle (0.5pt);
		\filldraw[Orchid] (1+1/6+1/6,-1) circle (2pt);
		\draw[->] (1+1/6+1/6,-1.5) node[anchor=north,black]{$w_{n_f}$} -- (1+1/6+1/6,-1.1);
		\filldraw[Orchid] (1+1/6+1/2+1/6,-1) circle (2pt) node[anchor=north,black]{};
		\filldraw[Orchid] (1+1/6+1+1/6,-1) circle (2pt) node[anchor=north,black]{};
		\filldraw[Black] (1/2+1/6+3/2+1/6+1/6,-1) circle (0.5pt);
		\filldraw[Black] (1/2+2/6+3/2+1/6+1/6,-1) circle (0.5pt);
		\filldraw[Black] (1/2+3/6+3/2+1/6+1/6,-1) circle (0.5pt);
		\filldraw[Orchid] (1+1/6+3/2+1/6+1/6,-1) circle (2pt);
		\draw[->] (1+1/6+3/2+1/6+1/6,-1.5) node[anchor=north,black]{$w_{n_f+n_g}$} -- (1+1/6+3/2+1/6+1/6,-1.1);
	\end{tikzpicture}
	\caption{}\label{fig:MergingVanKampen}
	\end{figure}
	From here, we observe that we obtain the same result by either merging $f$ and $g$ first and then applying \autoref{lem:LebesgueCoveringSC} or by first applying the lemma and then merging the function, proving that this relation does not affect $h$ in \autoref{Eq:DefinitionH}.
	
	The equivalence relation given by adding opposite maps ($\xleftrightarrow{d_{op}(i)}$) is a little more interesting, but it is easy to see that they cancel each other with $h$ defined as in \autoref{Eq:DefinitionH}.
	That is, consider the string $[f,\overline{f}]$, and apply \autoref{lem:LebesgueCoveringSC} to build $\gamma_{i,j}$ and $g_{i,j}$ as above. Thus we have that

\begin{align*}
	h[f,\overline{f}] = & \kappa_{1,0}[g_{1,L}] \star \kappa_{1,1}[g_{1,1}]\star \cdots \star \kappa_{1,R}[g_{1,R}] \\
	& \star \kappa_{1,R}[\overline{g_{1,R}}] \star \cdots \star \kappa_{1,1}[\overline{g_{1,1}}] \star \kappa_{1,L}[\overline{g_{1,L}}]\\
	= & \kappa_{1,0}[g_{1,L}] \star \kappa_{1,1}[g_{1,1}] \star \cdots \star \kappa_{1,R}[g_{1,R}] \\
	& \star \kappa_{1,R}[g_{1,R}]^{-1} \star \cdots \star \kappa_{1,1}[g_{1,1}]^{-1} \star \kappa_{1,L}[g_{1,L}]^{-1}.
\end{align*}
	It follows that a consecutive pair of opposite maps is mapped to the identity, and the equivalence classes are preserved.

	The most difficult case is delete-one-point semi-coarse homotopy equivalence ($\simeq_d$).
	Given their definition, it is enough to work with a one-string $[f]$.
	In the following, we first show that $h$ does not depend on the representative of delete-one-point semi-coarse homotopy equivalence.

	Let $f\xrightarrow{d(j)} g$.
	Consider the $z_k$ in the domain of $f$ (that is, $\Z_1$) as in \autoref{lem:LebesgueCoveringSC}.
	If $z_{k}< j < z_{k+1}$, there is no problem in the construction
	\autoref{lem:LebesgueCoveringSC} when we remove the element in $j$, and $h([f])=h([g])$.
	
	Now consider that $j=z_k$ for some $k$.
	Without loss of generality, suppose that $f(z_k)\in B$ (then $f(z_k)\notin A$ and $f(z_k-1)\in A$, by construction), and observe that $\{f(j-1), f(j+1))\}\in \mathcal{V}$, since $f \xrightarrow{d(j)} g$.
	
	In a first instance, we have two cases:
	\begin{enumerate}
	\item \label{Case:DefHom1} $f(z_k+1)\in B-A$,
	\item \label{Case:DefHom2} $f(z_k+1)\notin B-A$.
	\end{enumerate}
	If $f(z_k+1)\in B-A$, then by condition (1) of the definition of well-splitting, (\autoref{def:WellSplit}) there exists a short path $\gamma:\{0,1,2\}\to A\sqcup_{A\cap B}B$ such that $\gamma(0)=f(z_k-1)$ and $\gamma(2)=f(z_k+1)$.
	If, on the contrary, $f(z_k+1)\notin B-A$, then $z_{k-1}$ and $z_{k+1}$ are in the same factor defined in the construction in \autoref{lem:LebesgueCoveringSC}.
	
	In addition, for (1) in \autoref{def:WellSplit}, there exists $x'\in X$ such that
	\[ \{(x',f(z_k)), (x',f(z_k-1)), (x',f(z_k+1))\} \subset \mathcal{V}_{A\sqcup_{A\cap B} B}, \]
	and we can choose an $x'$ such that $x'\in A\cap B$ (see \autoref{fig:DeletingAPoint}, where we illustrate two ways to obtain $x'$ depending on the images of $z_k,z_k+1,z_k-1$.).
	
	\begin{figure}[h]
	\begin{tikzpicture}
	\filldraw[black] (0,0) circle (0pt) node[anchor=south]{$f(z_k)$};
	\filldraw[black] (0,-2) circle (0pt) node[anchor=north]{$f(z_k+1)$};
	\filldraw[black] (-3,-1) circle (0pt) node[anchor=east]{$f(z_k-1)$};
	\filldraw[black] (-1.35,-1) circle (0pt) node[anchor=west]{$x'$};

	\draw[Black,line width=0.5mm] (-3,-1) -- (0,0);
	\draw[Black,line width=0.5mm] (-3,-1) -- (0,-2);
	\draw[Red,line width=0.5mm] (-1.5,-1) -- (0,0);
	\draw[Red,line width=0.5mm] (-1.5,-1) -- (0,-2);
	\draw[Red,line width=0.5mm] (0,0) -- (0,-2);
	\draw[RoyalBlue,line width=0.5mm] (-3,-1) -- (-1.5,-1);

	\filldraw[Red] (0,0) circle (2pt);
	\filldraw[Red] (0,-2) circle (2pt);
	\filldraw[RoyalBlue] (-3,-1) circle (2pt);
	\filldraw[Orchid] (-1.5,-1) circle (2pt);
	
	\filldraw[black] (0+7,0) circle (0pt) node[anchor=south]{$f(z_k-1)$};
	\filldraw[black] (0+7,-2) circle (0pt) node[anchor=north]{$f(z_k+1)$};
	\filldraw[black] (-3+7,-1) circle (0pt) node[anchor=east]{$f(z_k)$};
	\filldraw[black] (-1.35+7,-1) circle (0pt) node[anchor=west]{$x'$};

	\draw[Black,line width=0.5mm] (-3+7,-1) -- (0+7,0);
	\draw[Black,line width=0.5mm] (-3+7,-1) -- (0+7,-2);
	\draw[RoyalBlue,line width=0.5mm] (-1.5+7,-1) -- (0+7,0);
	\draw[RoyalBlue,line width=0.5mm] (-1.5+7,-1) -- (0+7,-2);
	\draw[RoyalBlue,line width=0.5mm] (0+7,0) -- (0+7,-2);
	\draw[Red,line width=0.5mm] (-3+7,-1) -- (-1.5+7,-1);

	\filldraw[RoyalBlue] (0+7,0) circle (2pt);
	\filldraw[RoyalBlue] (0+7,-2) circle (2pt);
	\filldraw[Red] (-3+7,-1) circle (2pt);
	\filldraw[Orchid] (-1.5+7,-1) circle (2pt);
\end{tikzpicture}
\caption{Examples of adding $x'$. The left picture occurs when $f(z_k-1)\in A-B$ and $f(z_k+1)\in B-A$. The right picture occurs when $f(z_k-1)\in A-B$ and $f(z_k+1)\in A-B$.}
\label{fig:DeletingAPoint}
\end{figure}
	
	If $f(z_k+1)\in B-A$, then by condition (2) of well-split, we can select $\gamma(1)$ as $x'$, and therefore we can add $f(z_k)$ previous $f(z_k+1)$.
	We obtain what we want to because we have already proven that $h$ does not depend on the selection of $\gamma$.
	
	On the other hand, if $f(z_k+1)\notin B-A$, then we add $x'$ between $f(z_k-1)$ and $f(z_k+1)$ and split with respect that point.
	Now we can add $f(z_k)$ between $x'$ and $f(z_k+1)$ (see \autoref{fig:DeletingAPoint2}).
	We now obtain that $h$ is independent of the choice of $f$ and $g$ by also splitting the function at the point $z_{k+1}$, as desired.	

\begin{figure}[h!]
\begin{tikzpicture}
\foreach \x in {-3,...,-1}
	{\filldraw[black] (\x/6,0) circle (0.2pt);}
\foreach \x in {0,...,5}
	{\filldraw[RoyalBlue] (\x/2,0) circle (1.5pt);}
\foreach \x in {0,...,5}
	{\filldraw[Red] (\x/2+6/2,0) circle (1.5pt);}
\foreach \x in {0,...,3}
	{\filldraw[RoyalBlue] (\x/2+12/2,0) circle (1.5pt);}
\foreach \x in {0,...,2}
	{\filldraw[black] (\x/6+15/2+1/6,0) circle (0.2pt);}

\filldraw[black] (5/2,0) circle (0pt) node[anchor=north]{\footnotesize{$-1$}};	
\filldraw[black] (3,0) circle (0pt) node[anchor=north]{\footnotesize{$0$}};
\filldraw[black] (7/2,0) circle (0pt) node[anchor=north]{\footnotesize{$1$}};

\foreach \x in {-3,...,-1}
	{\filldraw[black] (\x/6,-1) circle (0.2pt);}
\foreach \x in {0,...,5}
	{\filldraw[RoyalBlue] (\x/2,-1) circle (1.5pt);}
\foreach \x in {0,...,4}
	{\filldraw[Red] (\x/2+6/2,-1) circle (1.5pt);}
\foreach \x in {0,...,3}
	{\filldraw[RoyalBlue] (\x/2+11/2,-1) circle (1.5pt);}
\foreach \x in {0,...,2}
	{\filldraw[black] (\x/6+14/2+1/6,-1) circle (0.2pt);}

\filldraw[black] (5/2,-1) circle (0pt) node[anchor=north]{\footnotesize{$-1$}};	
\filldraw[black] (3,-1) circle (0pt) node[anchor=north]{\footnotesize{$1$}};

\foreach \x in {-3,...,-1}
	{\filldraw[black] (\x/6,-2) circle (0.2pt);}
\foreach \x in {0,...,5}
	{\filldraw[RoyalBlue] (\x/2,-2) circle (1.5pt);}
\filldraw[Orchid] (6/2,-2) circle (1.5pt);
\foreach \x in {0,...,4}
	{\filldraw[Red] (\x/2+7/2,-2) circle (1.5pt);}
\foreach \x in {0,...,3}
	{\filldraw[RoyalBlue] (\x/2+12/2,-2) circle (1.5pt);}
\foreach \x in {0,...,2}
	{\filldraw[black] (\x/6+15/2+1/6,-2) circle (0.2pt);}

\filldraw[black] (5/2,-2) circle (0pt) node[anchor=north]{\footnotesize{$-1$}};	
\filldraw[black] (3,-2) circle (0pt) node[anchor=south]{\footnotesize{$x'$}};
\filldraw[black] (7/2,-2) circle (0pt) node[anchor=north]{\footnotesize{$1$}};

\foreach \x in {-3,...,-1}
	{\filldraw[black] (\x/6,-3) circle (0.2pt);}
\foreach \x in {0,...,5}
	{\filldraw[RoyalBlue] (\x/2,-3) circle (1.5pt);}
\filldraw[Orchid] (6/2,-3) circle (1.5pt);
\foreach \x in {0,...,5}
	{\filldraw[Red] (\x/2+7/2,-3) circle (1.5pt);}
\foreach \x in {0,...,3}
	{\filldraw[RoyalBlue] (\x/2+13/2,-3) circle (1.5pt);}
\foreach \x in {0,...,2}
	{\filldraw[black] (\x/6+16/2+1/6,-3) circle (0.2pt);}

\filldraw[black] (5/2,-3) circle (0pt) node[anchor=north]{\footnotesize{$-1$}};	
\filldraw[black] (3,-3) circle (0pt) node[anchor=south]{\footnotesize{$x'$}};
\filldraw[black] (7/2,-3) circle (0pt) node[anchor=north]{\footnotesize{$0$}};
\filldraw[black] (8/2,-3) circle (0pt) node[anchor=north]{\footnotesize{$1$}};
\end{tikzpicture}
\caption{Representation of the process; $f_i(z_{k+1}+1)\notin B-A$ and $z_k=0$.}
\label{fig:DeletingAPoint2}
\end{figure}
\end{proof}

\begin{teorema}
\label{theo:VankampenDisconnected}
Let $X$ be a semi-coarse space and $\mathcal{U}\subset 2^X$, and write $\mathcal{U}'$ to be the set of connected components of the elements of $\mathcal{U}$. Consider $A,B\subset X$ such that
\begin{enumerate}
\item $A,B$ well-split $X$,
\item for every component $C$ of $A$, $B$, and  $A\cap B$,  there exists a $U'\in\mathcal{U}'$ such that $C\cap U'\neq \varnothing$, and
\item for every $U'\in\mathcal{U}'$, the statement $(U'\subset A) \vee (U'\subset B)$ is true.
\end{enumerate}
Then the diagram
\begin{align*}
\xymatrix{
\pi_{\leq 1}(A\cap B, \mathcal{U}\cap A\cap B) \ar[r]^{i_A} \ar[d]^{i_B} &  \pi_{\leq 1}(A, \mathcal{U}\cap A) \ar[d]^{i_1}\\
\pi_{\leq 1}(B, \mathcal{U}\cap B) \ar[r]^{i_2} & \pi_{\leq 1}(X, \mathcal{U})
}
\end{align*}
is a pushout.
\end{teorema}

\begin{proof}
Divide $X$ in all its components and apply $X_\alpha$ and apply \autoref{theo:VanKampenConnected} to every $X_\alpha$ and $X_\alpha\cap \mathcal{U}'$, obtaining the result.
\end{proof}

We conclude this section with two consequences of this theorem.

\begin{corolario}
\label{Coro:VanKampenOneSet}
Let $X$ be a semi-coarse space and $\mathcal{U}=\{U\}$ for some $U\subset X$. Write $\mathcal{U}'$ the collection of all of the components of $U$. Consider $A,B\subset X$ such that
\begin{enumerate}
\item $A,B$ well-split $X$,
\item for every component $C$ of $A$, $B$, and  $A\cap B$, there exists a $U'\in\mathcal{U}'$ such that $C\cap U'\neq \varnothing$, and
\item for every $U'\in\mathcal{U}'$, the statement $(U'\subset A) \vee (U'\subset B)$ is true.
\end{enumerate}

Then the diagram
\begin{align*}
\xymatrix{
\pi_{\leq 1}(A\cap B, \mathcal{U}'\cap A\cap B) \ar[r]^{i_A} \ar[d]^{i_B} &  \pi_{\leq 1}(A, \mathcal{U}'\cap A) \ar[d]^{i_1}\\
\pi_{\leq 1}(B, \mathcal{U}'\cap B) \ar[r]^{i_2} & \pi_{\leq 1}(X, \mathcal{U})
}
\end{align*}
is a pushout.
\end{corolario}

\begin{corolario}
\label{Coro:VanKampenDiscrete}
Let $X$ be a semi-coarse space and $\mathcal{U}=\{\{x\} \mid x\in X\}$. Consider $A,B\subset X$ such that $A,B$ well-split $X$. Then the diagram
\begin{align*}
\xymatrix{
\pi_{\leq 1}(A\cap B, \mathcal{U}\cap A\cap B) \ar[r]^{i_A} \ar[d]^{i_B} &  \pi_{\leq 1}(A, \mathcal{U}\cap A) \ar[d]^{i_1}\\
\pi_{\leq 1}(B, \mathcal{U}\cap B) \ar[r]^{i_2} & \pi_{\leq 1}(X, \mathcal{U})
}
\end{align*}
is a pushout. In addition, the groupoid $\pi_{\leq 1}(X,\mathcal{U})$ contains all of the $\pi_1^{sc}(X,x)$.
\end{corolario}

\appendix

\section{Coarse Homotopy Fails in Semi-Coarse}
\label{Appex:CoarseHomotopy}

Many coarse invariants have been studied before (for example, homology theories \cite{Bunke_2020}, homotopy \cite{Mitchener_2020}, growth and amenability \cite{Roe_2003}), and a natural question is whether we may adapt them to the study of semi-coarse spaces. The problems with this kind of procedure are diverse, but in general, if we adapt the definition of a coarse invariant to semi-coarse spaces, then in some cases we do not recover the same information for coarse spaces, and in others we lose some important properties of the invariants. In this brief appendix, we examine one such example, the \emph{coarse homotopy} studied in  \cite{Mitchener_2020}. In addition, not all semi-coarse invariants are naturally useful when restricted to coarse spaces, As an example, recall that the semi-coarse homotopy groups studied in \cite{rieser2023semicoarse} are trivial in coarse spaces (\cite{rieser2023semicoarse}, Theorem 3.2.19).

In this section, $\mathbb{R}_+$ is endowed  with the coarse structure
\begin{align*}
\{A\subset \mathbb{R}_+\times\mathbb{R}_+\mid \exists r>0 \text{ such that } d(x,y)\leq r\ \forall (x,y)\in A\}.
\end{align*}
In addition, we recall the definition of a coarse map.

\begin{definicion}
\label{def:BoundedSetsProperMapsCoarseMaps}
Let $(X,\mathcal{V})$ and $(Y,\mathcal{W})$ be coarse spaces
\begin{itemize}
\item $B\subset X$ is \emph{bounded in $(X,\mathcal{V})$} if there exists $x\in X$ and $A\in \mathcal{V}$ such that
\begin{align*}
A[x]\coloneqq \{ y\in X \mid (x,y)\in A \} = B.
\end{align*}
\item A map $f:X\rightarrow Y$ is \emph{proper} if $f^{-1}(B)$ is bounded in $(X,\mathcal{V})$ for every $B$ bounded in $(Y,\mathcal{W})$.
\item A map $f:X\rightarrow Y$ bornologous and proper is called \emph{coarse}.
\end{itemize}
\end{definicion}

\begin{definicion}[$p$-Cylinder; \cite{Mitchener_2020}, 2.1]
\label{def:pCylinder}
Let $X$ be a coarse space, and let $p:X\rightarrow \mathbb{R}_+$ be a coarse map. Then we define the \emph{$p$-cylinder}
\begin{align*}
I_p X\coloneqq \{ (x,t)\in X\times \mathbb{R}_+ \mid t\leq p(x)+1 \}
\end{align*}
We have inclusions $i_0:X\rightarrow I_p X$ and $i_1:X\rightarrow I_p X$ defined by the formulas $i_0(x)=(x,0)$ and $i(x)=(x,p(x)+1)$, respectively. We also have the canonical projection $q:I_p X\rightarrow X$ defined by the formula $q(x,t)=x$.
\end{definicion}

\begin{definicion}[Coarse Homotopy; \cite{Mitchener_2020}, 2.2]
\label{def:CoarseHomotopy}
Let $X$ and $Y$ be coarse spaces. A \emph{coarse homotopy} is a coarse map $H:I_p X\rightarrow Y$ for some coarse map $p:X\rightarrow \mathbb{R}_+$.

We call coarse maps $f_0: X\rightarrow Y$ and $f_1:X\rightarrow Y$ \emph{coarsely homotopic} if there is a coarse map $H:I_p X\rightarrow Y$ such that $f_0 = H\circ i_0$ and $f_1 = H_1\circ i_1$.

This map $H$ is termed a \emph{coarse homotopy} between the maps $f_0$ and $f_1$. 
\end{definicion}

	We explore the following three modifications of these notions for semi-coarse spaces. Let
$(X,\mathcal{V})$ be a semi-coarse space.

	First, as in coarse spaces, say that a set $A \subset X$ in a semi-coarse space is bounded if there exists a controlled set $B \in \mathcal{V}$ and $x\in X$ such that $A=B[x]$.
	The next example shows that there aren't  any coarse maps to develop homotopy as maps from $Z_1 \to X$ when $X$ is coarse, while this is .

\begin{ejemplo}
	Let $X$ be a coarse space and consider $f:\mathbb{Z}_1\rightarrow X$ be a coarse map (a bornologous proper map.)
	Since $(k,k+1)$ is a controlled set in $\mathbb{Z}_1$, then we have that $\{f(0),f(5)\}$ is a bounded set in $X$, however $\{0,5\}$ is not in $\mathbb{Z}_1$.
	Thus, we don't have any proper map from $\mathbb{Z}_1$ to $X$.
	Therefore we don't have any $p$-Cylinder defined for $\Z_1$, and then we cannot define this kind of homotopy.
	(We can mimic this procedure replacing $\mathbb{Z}_1$ for many semi-coarse spaces which are not coarse spaces.)
\end{ejemplo}

	Second, considering bornologous maps instead of proper maps, then we obtain coarse spaces which are equivalent with this new definition but they are not homotopy equivalent with the original example.

\begin{ejemplo}
	Consider $\mathbb{Z}$ with the coarse structure induced by the metric.
	We observe that the map $p:\mathbb{Z}\rightarrow \mathbb{R}_+$ such that $p(x)=|x|$ is bornologous (actually it is a coarse map).
	Now we define $H:I_p(\Z)\rightarrow \Z$ such that 
	\[H(x,t)= \begin{cases}
	0 & 0\leq t \leq 0.1,\\
	\lfloor x-0.1 \rfloor & t>0.1,\ x\geq 0,\\
	\lceil -x+0.1 \rceil & t>0.1,\ x< 0;
	\end{cases}\]
	$H$ is a bornologous map and we can observe that it is not proper. Additionally, we observe that $H\circ i_0 (x,t)=0$ and $H\circ i_1 (x,t)=x$.
	Then we found a homotopy between $id_\Z$ and the constant map $0$;
	thus $\Z$ and $\{0\}$ are homotopic, which does not happen in the original homotopy.
\end{ejemplo}

	 In addition, in the same case, considering bornologous maps, in semi-coarse spaces, we lose transitivity of the coarse homotopies.

\begin{ejemplo}
Consider $X$ as the semi-coarse space induced by the graph
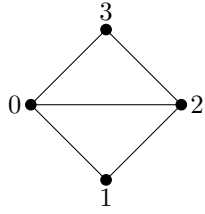
\begin{center}
\begin{tikzpicture}
\draw[black] (0,0) -- (1,1);
\draw[black] (0,0) -- (2,0);
\draw[black] (0,0) -- (1,-1);
\draw[black] (2,0) -- (1,1);
\draw[black] (2,0) -- (1,-1);

\filldraw[black] (0,0) node[anchor=east]{$0$} circle (2pt);
\filldraw[black] (1,1) node[anchor=south]{$3$} circle (2pt);
\filldraw[black] (2,0) node[anchor=west]{$2$} circle (2pt);
\filldraw[black] (1,-1) node[anchor=north]{$1$} circle (2pt);
\end{tikzpicture}
\end{center}
and the bornologous map (which is not proper) $p:\mathbb{Z}_1\rightarrow X$ such that $p(z)=0$. Then we obtain that $I_p(\mathbb{Z}_1)=\mathbb{Z}_1\times [0,1]$. Let's define the following bornologous functions:
\begin{align*}
u:\mathbb{Z}_1\rightarrow X & \text{ such that } u(3n)=0,\ u(3n+1)=2,\ u(3n+2)=3 \text{ for every }n\in\mathbb{Z}\\
m:\mathbb{Z}_1\rightarrow X & \text{ such that } u(3n)=0,\ u(3n+1)=u(3n+2)=2 \text{ for every }n\in\mathbb{Z}\\
d:\mathbb{Z}_1\rightarrow X & \text{ such that } u(3n)=0,\ u(3n+1)=2,\ u(3n+2)=1 \text{ for every }n\in\mathbb{Z}
\end{align*}
Observe that we have the following homotopies:
\begin{align*}
H_1(z,t)= u(z) \text{ if } 0\leq t <1,\ H_1(z,t)= m(z) \text{ if } t=1,\\
H_2(z,t)= d(z) \text{ if } 0\leq t <1,\ H_2(z,t)= m(z) \text{ if } t=1
\end{align*}
However, there is no $H$ such that $H\circ i_0=u$ and $H\circ i_1 = d$ because ${1,3}$ is not controlled in $X$.
\end{ejemplo}

	In the last case, we could try to replace $\mathbb{R}_+$ for a semi-coarse space in the definition of $I_p(X)$ or in the codomain of $H$.
	Observe that if $(X,\mathcal{V})$ is a connected coarse space, $(Y,\mathcal{W})$ is a semi-coarse space and $f:X\rightarrow Y$ is a bornologous map, then $f(X)$ is a connected (in the coarse sense) space.
	This is a problem because it might limit our bornologous maps. To mention a brief example, if you replace $\mathbb{R}_+$ by $\mathbb{R}_{+,r}$ you will observe that every bornologous function $p$ maps, for example, $\Z$ to $p(\Z)$ with coarse structure $2^{p(\Z)\times p(\Z)}$.
	Thus, the controlled sets in $I_p(\Z)$ satisfies that $\{ (z,t), (z',t') \}$ if and only if $\{(z,z')\}$ is controlled in $\Z$.
	Therefore, the image every bornologous maps $H:I_p(\Z)\to \Z_1$ necessary has to be contained in $\{k,k+1\}$ for some $k\in\Z$.
	It is easy to observe that the constant maps $f:\Z\to \Z_1$ such that $z\mapsto 1$, $g:\Z\to \Z_1$ such that $z\mapsto 2$ and $h:\Z\to \Z_1$ such that $z\mapsto 3$ satisfies that $f$ is coarse homotopic to $g$and $g$ are coarse homotpic to $h$, but $f$ is not coarse homotopic to $h$.

\section*{Acknowledgments}

	The author would like to thank Antonio Rieser, who is the supervisor of my PhD project, for the discussions about this paper and the possible directions this work could take, as well as his advice to include a brief discussion about the coarse homotopy groups.
	The author would also like to thank Noe Barcenas for his feedback on an event at the Casa Mexicana de Matemáticas in Oaxaca and another event in the CCM in Morelia.

\printnomenclature
\bibliographystyle{amsalpha}
\bibliography{all}


\end{document}